\documentclass[12pt]{article}
\usepackage{graphicx}
\usepackage{amsmath,amsthm,amssymb,enumerate}
\usepackage{euscript,mathrsfs}
\usepackage{dsfont}
\usepackage{xcolor}
\usepackage{empheq}
\usepackage[left=2cm,right=2cm,top=3.5cm,bottom=3.5cm]{geometry}
\usepackage{color}
\catcode`\@=11 \@addtoreset{equation}{section}

\catcode`\@=12

\allowdisplaybreaks

\newtheorem{Theorem}{Theorem}[section]
\newtheorem{Proposition}[Theorem]{Proposition}
\newtheorem{Lemma}[Theorem]{Lemma}
\newtheorem{Corollary}[Theorem]{Corollary}

\theoremstyle{definition}
\newtheorem{Definition}[Theorem]{Definition}

\newtheorem{Remark}[Theorem]{Remark}

\newcommand{\bTheorem}[1]{
\begin{Theorem} \label{T#1} }
\newcommand{\eT}{\end{Theorem}}

\newcommand{\bProposition}[1]{
\begin{Proposition} \label{P#1}}
\newcommand{\eP}{\end{Proposition}}

\newcommand{\bLemma}[1]{
\begin{Lemma} \label{L#1} }
\newcommand{\eL}{\end{Lemma}}
\newcommand{\bCorollary}[1]{
\begin{Corollary} \label{C#1} }
\newcommand{\eC}{\end{Corollary}}

\newcommand{\vrh}{\vr_h}
\newcommand{\vuh}{\vu_h}
\newcommand{\mh}{\vc{m}_h}
\newcommand{\Dh}{\frac{{\rm d}}{{\rm d}t} }
\newcommand{\bRemark}[1]{
\begin{Remark} \label{R#1} }
\newcommand{\eR}{\end{Remark}}

\newcommand{\bDefinition}[1]{
\begin{Definition} \label{D#1} }
\newcommand{\eD}{\end{Definition}}

\newcommand{\Del}{\Delta_x}

\newcommand{\intSh}[1] {\int_{\sigma} #1 {\rm d}S_h }

\newcommand{\bfphi}{\boldsymbol{\varphi}}
\newcommand{\bfPhi}{\boldsymbol{\Phi}}

\newcommand{\bFormula}[1]{
\begin{equation} \label{#1}}
\newcommand{\eF}{\end{equation}}

\newcommand{\grid}{\mathcal{T}_h}

\newcommand{\Ov}[1]{\overline{#1}}

\newcommand{\aleq}{\stackrel{<}{\sim}}
\newcommand{\ageq}{\stackrel{>}{\sim}}

\newcommand{\vr}{\varrho}

\newcommand{\vt}{\vartheta}
\newcommand{\vu}{\vc{u}}
\newcommand{\vm}{\vc{m}}

\newcommand{\vn}{\vc{n}}

\newcommand{\vc}[1]{{\bf #1}}

\newcommand{\Div}{{\rm div}_x}
\newcommand{\Grad}{\nabla_x}
\newcommand{\Dt}{\frac{\rm d}{{\rm d}t}}

\newcommand{\dx}{\,{\rm d} {x}}

\newcommand{\dt}{\,{\rm d} t }

\newcommand{\intO}[1]{\int_{\Omega} #1 \ \dx}
\newcommand{\intOh}[1]{\int_{\Omega_h} #1 \ \dx}

\def\softd{{\leavevmode\setbox1=\hbox{d}%
          \hbox to 1.05\wd1{d\kern-0.4ex{\char039}\hss}}}

\definecolor{Cgrey}{rgb}{0.85,0.85,0.85}
\definecolor{Cblue}{rgb}{0.50,0.85,0.85}
\definecolor{Cred}{rgb}{1,0,0}
\definecolor{fancy}{rgb}{0.10,0.85,0.10}

\newcommand\Cbox[2]{%
    \newbox\contentbox%
    \newbox\bkgdbox%
    \setbox\contentbox\hbox to \hsize{%
        \vtop{
            \kern\columnsep
            \hbox to \hsize{%
                \kern\columnsep%
                \advance\hsize by -2\columnsep%
                \setlength{\textwidth}{\hsize}%
                \vbox{
                    \parskip=\baselineskip
                    \parindent=0bp
                    #2
                }%
                \kern\columnsep%
            }%
            \kern\columnsep%
        }%
    }%
    \setbox\bkgdbox\vbox{
        \color{#1}
        \hrule width  \wd\contentbox %
               height \ht\contentbox %
               depth  \dp\contentbox
        \color{black}
    }%
    \wd\bkgdbox=0bp%
    \vbox{\hbox to \hsize{\box\bkgdbox\box\contentbox}}%
    \vskip\baselineskip%
}

\date{}

\begin{document}


\title{A finite volume scheme for the Euler system \\inspired by the two velocities approach}

\author{Eduard Feireisl
\thanks{The research of E.F. and H.M.~leading to these results has received funding from
the Czech Sciences Foundation (GA\v CR), Grant Agreement
18--05974S. The Institute of Mathematics of the Academy of Sciences of
the Czech Republic is supported by RVO:67985840.} \and M\' aria Luk\' a\v cov\' a-Medvi\softd ov\'a
\thanks{The reserach of M.L. was supported by the German Science Foundation under the Collaborative Research Centers TRR~146 and TRR~165.}
\and Hana Mizerov\' a \footnotemark[1] $^\mathsection$
}

\date{\today}

\maketitle

\bigskip

\centerline{$^*$ Institute of Mathematics of the Academy of Sciences of the Czech Republic}

\centerline{\v Zitn\' a 25, CZ-115 67 Praha 1, Czech Republic}

\centerline{feireisl@math.cas.cz}
\centerline{mizerova@math.cas.cz}

\bigskip

\centerline{$^\dagger$ Institute of Mathematics, Johannes Gutenberg-University Mainz}

\centerline{Staudingerweg 9, 55128 Mainz, Germany}

\centerline{lukacova@uni-mainz.de}

\bigskip
\centerline{$^\mathsection$ Department of Mathematical Analysis and Numerical Mathematics}
\centerline{Faculty of Mathematics, Physics and Informatics, Comenius University in Bratislava}

\centerline{Mlynsk\' a dolina, 842 48 Bratislava, Slovakia}

\begin{abstract}

We propose a  new finite volume scheme for the Euler system of gas dynamics motivated by the model proposed by H.~Brenner. Numerical
viscosity imposed through upwinding acts on the velocity field rather than on the convected quantities. The resulting numerical method
enjoys the crucial properties of the Euler system, in particular positivity of the approximate density and pressure and the minimal entropy principle. In addition, the approximate solutions generate a dissipative measure--valued solutions of the limit system.
 In particular, the numerical solutions converge to the smooth solution of the system as long as the latter exists.

\end{abstract}


\tableofcontents

\section{Introduction}
\label{b}

In 2005, H.~Brenner~\cite{BREN} proposed a new approach to dynamics of \emph{viscous and heat conducting} fluids based on two velocity fields distinguishing the bulk
mass transport from the purely microscopic motion. Brenner's approach has been subjected to thorough criticism by
\" Ottinger et al.~\cite{OeStLi}, where its incompatibility with certain physical principles is shown.  Nevertheless, some computational simulations have been performed by Greenschields and Reese~\cite{GrRe}, Bardow and \" Ottinger~\cite{BaOett}, Guo and Xu~\cite{GuoXu}
showing suitability of the model in specific situations.
More recently, Guermond and Popov~\cite{GuePop}
rediscovered the model pointing out its striking similarity with certain numerical methods based on the finite volume approximation
of the \emph{inviscid} fluids. In particular, unlike the conventional and well accepted Navier--Stokes--Fourier system,
Brenner's model reflects the basic properties of the complete Euler system in the asymptotic limit of vanishing transport coefficients.

Inspired by these observations, we propose a new finite volume scheme for the complete Euler system based on Brenner's ideas.
In particular, the new scheme enjoys the following properties:
\begin{itemize}
\item {\bf Positivity of the discrete density and temperature}\\ The approximate density and temperature remain strictly positive on any finite time interval.
\item {\bf Entropy stability}\\ The discrete entropy inequality  in the sense of Tadmor is satisfied, see~\cite{Tad87,Tad03}.
\item {\bf Minimum entropy principle} \\ The entropy attains its minimum at the initial time, cf. \cite{Tad86,GuePop}.
\item {\bf Weak BV estimates}\\ We control suitable weak BV norms of the discrete density, temperature and velocity.
\end{itemize}
In comparison with the conventional convergence results based  on unrealistic hypothesis on uniform boundedness of all physical quantities our scheme produces convergent solutions
as long as the gas remains in its non--degenerate regime, cf. Section~\ref{L}.

\subsection{Complete Euler system}

The complete Euler system describes the time evolution of the \emph{standard} physical fields: the mass density
$\vr = \vr(t,x)$, the macroscopic velocity $\vu = \vu(t,x)$, and the (absolute) temperature $\vt = \vt(t,x)$ of a perfect
compressible fluid,
\begin{align*}
\partial_t \vr + \Div (\vr \vu) &= 0,\\
\partial_t (\vr \vu) + \Div (\vr \vu \otimes \vu) + \Grad p &= 0,\\
\partial_t \left( \frac{1}{2} \vr |\vu|^2 + \vr e \right) +
\Div \left[ \left( \frac{1}{2} \vr |\vu|^2 + \vr e + p \right) \vu \right] &= 0.
\end{align*}
For the sake of simplicity, we consider the standard polytropic EOS with the Boyle--Marriot pressure law,
\[
p = (\gamma - 1) \vr e = \vr \vt, \ e = c_v \vt, \ c_v = \frac{1}{\gamma - 1}.
\]
Accordingly, the \emph{physical} entropy reads
\[
s(\vr, \vt) = \log \left( \frac{\vt^{c_v}}{\vr} \right)
\]
with the associated entropy inequality,
\[
\partial_t (\vr s) + \Div (\vr s \vu) \geq 0.
\]
Note that the same inequality is automatically satisfied by any ``renormalized'' \emph{mathematical} entropy $s_\chi$
\[
s_\chi = \chi \left( \log \left( \frac{\vt^{c_v}}{\vr} \right) \right),
\]
where $\chi$ is a non--decreasing concave function.

Numerical schemes are based on the \emph{conservative} variables:
the density $\vr$, the momentum $\vc{m} = \vr \vu$, and the total energy
\[
E = \frac{1}{2} \vr |\vu|^2 + \vr e.
\]
Accordingly, the Euler system takes the form
\begin{align}
\partial_t \vr + \Div \vc{m} = 0,\label{E7}\\
\partial_t \vc{m} + \Div \left( \frac{\vc{m} \otimes \vc{m} }{\vr} \right) + \Grad p = 0,\label{E8}\\
\partial_t E + \Div \left[ (E + p) \frac{\vc{m}}{\vr} \right] = 0, \label{E9}
\end{align}
where
\[
p = (\gamma - 1) \left( E - \frac{1}{2} \frac{|\vc{m}|^2}{\vr} \right).
\]
In the conservative framework, positivity of the density as well as of the pressure becomes an issue, in which the associated entropy
balance
\[
\partial_t (\vr s_\chi) + \Div \left( s_{\chi} \vc{m} \right) \geq 0
\]
plays a crucial role.

\subsection{Brenner's model}

Brenner's approach to modelling real \emph{viscous} and \emph{heat conducting} fluids postulates two velocities
$\vu$ and $\vc{v}$ interrelated through
\[
\vc{v} = \vu - K \nabla_x\log(\vr).
\]
For the Newtonian viscous stress
\[
\mathbb{S}(\Grad \vu) = \eta_1 \left( \Grad \vu + \Grad^t \vu - \frac{2}{3} \Div \vu \mathbb{I} \right) +
\eta_2 \Div \vu \mathbb{I}
\]
and the Fourier heat flux
\[
\vc{q} = -\kappa \Grad \vt
\]
the Brenner model reads
\begin{align}
\partial_t \vr + \Div (\vr \vc{v} ) &= 0, \label{E11}\\
\partial_t (\vr \vu) + \Div (\vr \vu \otimes \vc{v}) + \Grad p &= \Div \mathbb{S}(\Grad \vu),  \label{E12}\\
\partial_t \left( \frac{1}{2} \vr |\vu|^2 + \vr e \right) +
\Div \left[ \left( \frac{1}{2} \vr |\vu|^2 + \vr e + p \right) \vc{v} \right]
+ \Div \vc{q} &= \Div \left( \mathbb{S}(\Grad \vu) \cdot \vu \right),\label{E13}
\end{align}
see Brenner \cite{BREN2,BREN,BREN1}.
Moreover, if $K$ is related to the heat conductivity coefficient $\kappa$ through
\[
K = \frac{\kappa}{\vr c_v},
\]
then the associated entropy balance takes the form
\begin{equation}
\begin{aligned}\label{entropy_cont}
\partial_t (\vr s_\chi) &+ \Div (\vr s_{\chi} \vc{v} ) - \Div \left( \frac{\kappa}{c_v} \Grad s_\chi \right) \\
&= \frac{\chi'(s)}{\vt} \mathbb{S}(\Grad \vu) : \Grad \vu +
\kappa \chi'(s) |\Grad \log(\vt) |^2 + \chi'(s) \frac{\kappa}{c_v} |\Grad \log(\vr)|^2
- \chi''(s) \frac{\kappa}{c_v} |\Grad s |^2,
\end{aligned}
\end{equation}
see Guermond and Popov \cite{GuePop} and \cite[Section 4.1]{BreFei17a}.

As observed by Guermond and Popov \cite{GuePop}, for the ansatz
\[
\mathbb{S}(\Grad \vu) = h  \lambda \vr \Grad \vu + h^\alpha \Grad \vu,\
\kappa = c_v \vr K = c_v h \vr  \lambda, \  \lambda \geq 0,\]
the system (\ref{E11}--\ref{E13}) rewrites in the conservative variables as
\begin{align}
\partial_t \vr + \Div (\vr \vu) &= h \Div ( \lambda \Grad \vr),\label{b4}\\
\partial_t \vc{m} + \Div (\vc{m} \otimes \vu ) + \Grad p &=
h \Div \left(  \lambda \Grad \vc{m} \right) + h^\alpha \Del \vu, \label{b5}\\
\partial_t E + \Div (E \vu + p \vu) &= h \Div ( \lambda \Grad E) + h^\alpha \Div (\Grad \vu \cdot \vu ).\label{b6}
\end{align} 
This form, without the $h^{\alpha}$--dependent terms, is strongly reminiscent of some numerical schemes for the
complete (inviscid) Euler system based on the finite volume method like the Lax--Friedrichs scheme.

\subsection{Finite volume scheme}

Motivated by Guermond and Popov \cite{GuePop} we propose a finite volume scheme for the complete Euler system
based on (\ref{b4}--\ref{b6}). Although written exclusively in the conservative variables, the scheme relies on convective terms
expressed in terms of the velocity $\vu$ rather than the momentum $\vc{m}$. This allows to minimize the effect of the viscous perturbations
- a potential source of deviation from the target Euler system  for inviscid flows.
Indeed the scheme preserves all the basic properties of the continuous system,
in particular, it is entropy stable. Moreover,
the positivity of the density and pressure as well as the minimum entropy principle hold.

We then examine the properties of the associated semi--discrete dynamical system. We show  that it generates in the asymptotic limit a dissipative measure--valued (DMV) solution of the complete Euler system introduced in
\cite{BreFei17,BreFei17a}, see also \cite{FLM18} for the convergence of the Lax--Friedrichs method. Moreover, employing the (DMV)--strong uniqueness principle, we will obtain strong (pointwise) convergence to the unique classical solution as long as the latter exists.
In contrast with the standard entropy stable finite volume methods, where convergence analysis is based on rather unrealistic {\it a priori} hypotheses of uniform boundedness of numerical solutions, cf. Fjordholm, Mishra, K{\" a}ppeli, Tadmor \cite{FjKaMiTa, FjMiTa2, FjMiTa1,Tad03}, the convergence for the present scheme is almost unconditional, requiring only a technical hypothesis of boundedness of the numerical temperature and the absence of vacuum.

The paper is organized as follows. Section \ref{n} contains necessary preliminaries including the geometric properties of the mesh and the basic notation used in finite volume methods. Then we introduce the numerical method and the associated semi--discrete dynamical system.
In Section~\ref{S}, we show that the scheme is entropy  stable. In Section~\ref{AS},
we study stability of the semi--discrete scheme deriving all necessary {\it a priori} bounds.
Consistency of the scheme, based on a careful analysis of the error terms, is discussed in Section~\ref{C}. Finally, we perform the limit of vanishing numerical step in Section~\ref{L}.

\section{Numerical scheme}
\label{n}

We introduce the basic notation, function spaces, and, finally, the numerical scheme.

\subsection{Preliminaries}

We suppose the physical space to be a polyhedral domain $\Omega_h \subset R^N$, $N=1,2,3$, that is decomposed into compact elements
\[
 \Ov{\Omega}_h = \bigcup_{K \in \grid} K.
\]
The elements $K$ are sharing either a common face, edge, or vortex. The mesh $\grid$ satisfies the standard regularity assumptions, cf.~\cite{ciarlet,EyGaHe}. The set of all faces is denoted by $\Sigma,$ while $\Sigma_{int}=\Sigma\backslash\partial\Omega_h$ stands for the set of all interior faces. Each face
is associated with a normal vector $\vc{n}$.
In what follows, we shall suppose
\[
|K|_N \approx h^N, \ |\sigma|_{N-1} \approx h^{N-1} \ \mbox{for any}\ K \in \grid, \ \sigma \in \Sigma.
\]
The symbol $Q_h$ denotes the set of functions constant on each element $K$.
For a piecewise (elementwise) continuous function $v$ we define
\[
v^{\rm out}(x) = \lim_{\delta \to 0+} v(x + \delta \vc{n}),\
v^{\rm in}(x) = \lim_{\delta \to 0+} v(x - \delta \vc{n}),\
\Ov{v}(x) = \frac{v^{\rm in}(x) + v^{\rm out}(x) }{2},\
[[ v ]] = v^{\rm out}(x) - v^{\rm in}(x)
\]
whenever $x \in \sigma \in \Sigma_{int}$. We recall the product rule
\[
[[ u v ]] = \Ov{u} [[v]] + [[u]] \Ov{v}.
\]

For $\Phi \in L^1(\Omega_h)$ we define the projection
\[
\Pi_h [\Phi] = \sum_{K \in \grid} 1_{K} \frac{1}{|K|} \int_K \Phi \dx \in Q_h(\Omega_h).
\]
If $\Phi \in C^1(\Ov{\Omega_h})$ we have
\begin{equation} \label{n4c}
\Big| \ [[ \Pi_h [\Phi] ]]  \ \Big|_{\sigma} \aleq h \| \Phi \|_{C^1},\
\left| \Phi - \Ov{\Pi_h [\Phi] } \right|_{\sigma} \aleq h \| \Phi \|_{C^1}
\  \mbox{for any}\ x \in \sigma \in \Sigma_{int}.
\end{equation}
Here and hereafter the symbol $A \aleq B$ means $A \leq cB$ for a generic positive constant $c$ independent of $h.$
If $\Phi \in C^2(\Ov{\Omega_h})$ and $\grid$ consists of uniform rectangular/cubic elements, then we moreover have
\begin{equation} \label{n4c2}
\frac{1}{|\sigma|} \intSh{\left| \Phi - \Ov{\Pi_h [\Phi] } \right|}  \aleq h^2 \| \Phi \|_{C^2}
\quad  \mbox{for any} \ \sigma \in \Sigma_{int}.
\end{equation}
Indeed, any $C^2$ function can be approximated by the piecewise linear Rannacher--Turek elements \cite{RT} (an analogue of the 
Crouzeix--Raviart elements on rectangles) with the error of $\mathcal{ O}(h^2).$  Thus, it is enough to show (\ref{n4c2}) for the non--conforming piecewise linear Rannacher--Turek elements. Taking into account their continuity in the center of cell interfaces and the definition of projection $\Pi_h,$ we only need to show
$$
\left| \Phi(S_\sigma) - \dfrac{\left(\Phi(S_K) +  \Phi(S_L)\right)}{2}\right| \aleq h^2,
$$ where $S_\sigma$ denotes the center of gravity of $\sigma$,
$S_K$ and $S_L$ the centers of gravity of two neighbouring elements $K$ and $L$ sharing the common face $\sigma.$ The latter follows directly from
the Taylor expansion. 

We further recall the negative  $L^p$--estimates \cite{ciarlet}
\begin{equation} \label{n4b}
\| v \|_{L^p(\Omega_h)} \aleq h^{ N\frac{1 - p}{p} } \| v \|_{L^1(\Omega_h)} \ \mbox{for any}\ 1 \leq p \leq \infty,\
 \mbox{with}\ \frac{N(1 - p)}{p} = - N \ \mbox{if} \ p = \infty,
\end{equation}
and the trace inequality
\begin{align} \label{trace}
\| v \|_{L^p(\partial K)} \aleq h^{ -\frac{1}{p} } \| v \|_{L^p(K)} \ \mbox{for any}\ 1 \leq p \leq \infty,
\end{align}
for any $v \in Q_h (\Omega_h)$.
Moreover, we have a discrete version of the Sobolev embedding theorem, see Chainais--Hillairet, Droniou \cite[Lemma 6.1]{ChHD},
\begin{equation} \label{n4d}
\| v \|_{L^6(\Omega_h)} \aleq \| v \|_{L^2(\Omega_h)} +  \left(\sum_{\sigma \in \Sigma_{int}} \intSh{
\frac{[[ v ]]^2}{h} }\right)^{1/2} \ \mbox{for any}\ v \in Q_h(\Omega_h),\ N=1,2,3.
\end{equation}

Given a velocity $\vu \in Q_h(\Omega_h; R^N)$ and $r \in Q_h(\Omega_h)$, we define on each face  $\sigma \in \Sigma_{int}$ an \emph{upwind} of $r$ by $\vu$ as
\begin{align}\label{Up}
Up [r, \vu] = \Ov{r} \ \Ov{\vu} \cdot \vc{n} - \frac{1}{2} |\Ov{\vu} \cdot \vc{n}| [[r]] =
r^{\rm in} [\Ov{\vu} \cdot \vc{n}]^+ + r^{\rm out} [\Ov{\vu} \cdot \vc{n}]^-.
\end{align}
Finally,  we set
\begin{align}\label{up_down}
r^{\rm up} = \left\{ \begin{array}{ll} r^{\rm in} & \mbox{if} \ \Ov{\vu} \cdot \vc{n} \geq 0\\ \\
r^{\rm out} & \mbox{if} \ \Ov{\vu} \cdot \vc{n} < 0,
\end{array} \right., \qquad
r^{\rm down} = \left\{ \begin{array}{ll} r^{\rm out} & \mbox{if} \ \Ov{\vu} \cdot \vc{n} \geq 0\\ \\
r^{\rm in} & \mbox{if} \ \Ov{\vu} \cdot \vc{n} < 0,
\end{array} \right.,
\end{align}
and
\begin{align}\label{tilde_jump}
\widetilde{[[r]]}=r^{\rm up}-r^{\rm down}= -[[r]] \,\textnormal{sgn}(\Ov{\vu}\cdot\vn).
\end{align}

\subsection{Approximation scheme}

In order to properly define  the numerical scheme, the boundary conditions must be specified. Here, we adopt the no--flux boundary condition:
\begin{align*}
\vuh\cdot\vn=0, \ \mbox{ for any } \sigma\in \partial\Omega_h,
\end{align*}
and $\vrh,$ $p_h$ are extrapolated, i.e. $\partial \vrh /\partial \vn=0=\partial p_h/\partial \vn,$ $\vn$ is an outer normal to $\partial\Omega_h.$
 We consider the numerical flux function in the form
\begin{align}\label{num_flux}
F_h(r_h,\vu_h)={Up}[r_h, \vuh] - \mu_h [[ r_h ]],
\end{align}
where $\mu_h\geq 0$ and ${Up}[r_h, \vuh]$ is given by \eqref{Up}.
The quantities
 $\vrh \in Q_h(\Omega_h)$, $\vc{m}_h \in Q_h(\Omega_h; R^N)$, and $E_h \in Q_h(\Omega_h)$
at the time level $t$ are given by the following system of equations:
\begin{itemize}
\item {\bf Continuity equation}
\begin{equation} \label{n2}
\intOh{ D_t \vrh \Phi } - \sum_{ \sigma \in \Sigma_{int} } \intSh{  F_h(\vrh,\vuh)
[[\Phi]]  } = 0 \ \mbox{for any}\ \Phi \in Q_h (\Omega_h),
\end{equation}
where
\[
\vc{u}_h = \frac{\vc{m}_h}{\vr_h}.
\]
\item {\bf Momentum equation}
\begin{equation} \label{n3}
\begin{split}
\intOh{ D_t \mh \cdot \bfPhi } &- \sum_{ \sigma \in \Sigma_{int} } \intSh{ {\bf F}_h(\mh,\vuh)
\cdot [[\bfPhi]]  }- \sum_{ \sigma \in \Sigma_{int} } \intSh{ \Ov{p_h} \vc{n} \cdot [[ \bfPhi ]] } \\
&= - h^{\alpha - 1} \sum_{ \sigma \in \Sigma_{int} } \intSh{ [[ \vu_h ]] \cdot [[ \bfPhi ]] } \ \mbox{for all}\
\bfPhi \in Q_{h} (\Omega_h, R^N),
\end{split}
\end{equation}
where
\[
p_h = (\gamma - 1) \left( E_h - \frac{1}{2} \frac{|\vc{m}_h|^2}{\vr_h} \right).
\]

\item {\bf Energy equation}
\begin{equation} \label{n4}
\begin{split}
\intOh{ D_t E_h \Phi } &- \sum_{ \sigma \in \Sigma_{int} } \intSh{  F_h(E_h,\vu_h)
[[\Phi]]  }\\
&- \sum_{ \sigma \in \Sigma_{int} } \intSh{ \Ov{ p_h } [[ \Phi \vu_h ]] \cdot \vc{n} }
+ \sum_{ \sigma \in \Sigma_{int} } \intSh{ \Ov{ p_h \Phi} [[\vu_h]] \cdot \vc{n} } \\
&= - h^{\alpha - 1} \sum_{ \sigma \in \Sigma_{int} } \intSh{ [[ \vu_h ]] \cdot \Ov{\vu_h}   [[ \Phi ]] } \ \mbox{for all}\
\Phi \in Q_{h} (\Omega_h).
\end{split}
\end{equation}
\end{itemize}
Note that our upwinding  ${Up}[r_h, \vu_h], $  $r_h = \vr_h, \vc{m}_h, E_h$,  is
based only on the sign of the normal component of velocity, instead of the sign of the eigenvalues as
in the standard flux--vector splitting schemes. In addition, numerical diffusion term $-\mu_h [[ r_h ]]  $ is added to the numerical flux function.
The parameter $\mu_h \geq 0$ is typically of the following form
\[
\mu_h = h M (h, \Ov{\vr_h}, \Ov{\mh}, \Ov{E_h}),
\]
where $M$ is a continuous function.
Unlike the convective terms, the pressure terms are appropriately  averaged, cf.~(\ref{n3}), (\ref{n4}). We should note that the terms on the right--hand side of (\ref{n3}), (\ref{n4}) can be interpreted as the interior penalty terms for the velocity $\vu_h$ that are typically used in the discontinuous Galerkin approach.

In the purely discrete version of (\ref{n2}--\ref{n4}), the operator $D_t$ stands for
\[
 D_t r_h = \frac{r_h (t) - r_h(t - \Delta t)}{\Delta t},
\]
where $\Delta t > 0$ is the time step. In the semi--discrete setting considered in this paper, the functions
$[\vrh, \mh, E_h]$ are continuous functions of the time $t \in [0,T]$, and $D_t$ is interpreted as the standard differential operator,
\[
D_t = \Dt.
\]

\begin{Remark} \label{nR1}

By virtue of the \emph{product rule},
the integral
\[
h^{\alpha - 1} \sum_{ \sigma \in \Sigma_{int} } \intSh{ [[ \vu_h ]] \cdot \Ov{\vu_h}   [[ \Phi ]] } = \frac{h^{\alpha - 1}}{2} \sum_{ \sigma \in \Sigma_{int} } \intSh{ [[ \vu_h^2 ]]    [[ \Phi ]] }
\]
may be replaced by a more convenient expression
\[
h^{\alpha - 1} \sum_{ \sigma \in \Sigma_{int} } \intSh{ [[ \vu_h ]] \cdot [[ \Phi \vuh ]] } -
h^{\alpha - 1} \sum_{ \sigma \in \Sigma_{int} } \intSh{ [[ \vu_h ]]^2 \Ov{\Phi} }.
\]

\end{Remark}

\begin{Remark} \label{nR1a}

We point out that
\begin{align}\label{r2.2}
\sum_{ \sigma \in \Sigma_{int} } \intSh{ \Ov{ p_h } [[ \Phi \vu_h ]] \cdot \vc{n} }
- \sum_{ \sigma \in \Sigma_{int} } \intSh{ \Ov{ p_h \Phi} [[\vu_h]] \cdot \vc{n} } \ne
\sum_{ \sigma \in \Sigma_{int} } \intSh{ \Ov{ p_h } \ \Ov{ \vu_h }\cdot \vc{n} [[ \Phi ]] }
\end{align}
as one might expect.
Indeed, the left--hand side of \eqref{r2.2} equals to
\begin{align}\label{pressure}
\sum_{ \sigma \in \Sigma_{int} } \intSh{ \Ov{ p_h } \ \Ov{ \vu_h }\cdot \vc{n} [[ \Phi ]] } -\frac{1}{4}\sum_{ \sigma \in \Sigma_{int} } \intSh{[[p_h]] [[\vu_h]]\cdot \vn [[ \Phi ]]  }.
\end{align}
\end{Remark}

This paper is devoted to the semi--discrete version, where $[\vrh, \mh, E_h]$ are continuous
functions of time and
the approximate scheme (\ref{n2}--\ref{n4}) may be interpreted as a finite system of ODEs.
It follows from the standard ODE theory that for a given initial state
\[
\begin{split}
\vr_h(0) &= \vr_{0,h} \in Q_h(\Omega_h), \ \vr_{0,h} > 0,\
 \vm(0) = \vm_{0, h} \in Q_h(\Omega_h;R^N),\ E_h(0) = E_{0,h} \in Q_h(\Omega_h),\\
E_{0,h} &- \frac{1}{2} \frac{|\vc{m}_{0,h}|^2}{\vr_{0,h}} > 0,
\end{split}
\]
the semi--discrete system (\ref{n2}--\ref{n4}) admits a unique solution $[\vrh, \mh, E_h]$ defined on a maximal time interval $[0, T_{\rm max}),$ where
\begin{equation} \label{nn5}
\vrh (t) > 0, \ p_h(t) = (\gamma - 1) \left( E_h(t) - \frac{1}{2} \frac{|\vc{m}_h(t) |^2}{\vrh(t)} \right)
> 0 \ \mbox{for all}\ t \in [0,T_{\rm max}).
\end{equation}
In particular, the absolute temperature $\vt_h$ can be defined,
\[
\vt_h (t) = \frac{ p_h(t) }{\vr_h (t) } = \frac{ \gamma - 1}{\vr_h (t)} \left( E_h(t) - \frac{1}{2} \frac{|\vc{m}_h(t) |^2}{\vrh(t)} \right).
\]

As we show in Section \ref{AS}, the system (\ref{n2}--\ref{n4}) admits sufficiently strong {\it a priori} bounds that will guarantee
{\bf (i)} $T_{\rm max} = \infty,$ {\bf (ii)} validity of (\ref{nn5}) for any $ t\geq 0.$

\section{Entropy balance}
\label{S}

We derive a discrete analogue of the entropy balance \eqref{entropy_cont} associated to the semi--discrete system (\ref{n2}--\ref{n4}).

\subsection{Renormalization}

The process of renormalization requires multiplying the discrete equations by nonlinear functions of the unknowns.

\subsubsection{Continuity equation}

Multiplying the continuity equation (\ref{b4}) by $b'(\vr)$ we deduce its renormalized form
\[
\partial_t b(\vr) + \Div (b(\vr) \vu) + \Big( b'(\vr) \vr - b(\vr) \Big) \Div \vu =
h \Div (\lambda \Grad b(\vr)) - \lambda b''(\vr) |\Grad \vr|^2.
\]
Its discrete analogue (\ref{n2}) gives rise to
\begin{equation} \label{r1}
\begin{split}
&\intOh{ \Dt b(\vr_h) \Phi } - \sum_{ \sigma \in \Sigma_{int} } \intSh{ {Up}[b(\vr_h), \vu_h] [[\Phi]] }
+ \sum_{ \sigma \in \Sigma_{int} } \intSh{ \Ov{\vuh} \cdot \vc{n} \left[\left[ \Big( b(\vr_h) - b'(\vr_h) \vr_h \Big) \Phi \right]\right]}
\\ &= - \sum_{ \sigma \in \Sigma_{int} } \intSh{ \mu_h [[\vrh]] \ [[ b'(\vrh) \Phi ]] }
- \sum_{ \sigma \in \Sigma_{int} } \intSh{
 \Phi^{\rm down}
\Big( \widetilde{[[ b(\vrh) ]]} - b'( \vrh^{\rm down} ) \widetilde{[[ \vrh ]]} \Big)  | \Ov{\vuh} \cdot \vc{n} | }
,
\end{split}
\end{equation}
for any  $\Phi \in  Q_h(\Omega_h)$, see \cite[Section 4.1]{FeKaNo}. Here $r^{\rm down}$ and $\widetilde{[[r_h]]}$ are given by \eqref{up_down} and \eqref{tilde_jump}, respectively.

\subsubsection{Transport equation}

Under the assumption that $\vr$ satisfies (\ref{b4}), we consider a field $b$ satisfying
\[
\partial_t (\vr b) + \Div (\vr b \vu) = F.
\]
Multiplying the equation by $\chi'(b)$ we obtain
\[
\partial_t (\vr \chi(b)) + \Div (\vr \chi(b) \vu) = F \chi'(b) + \Div ( h  \mu  \Grad \vr) \left( \chi(b) - b \chi'(b) \right).
\]
The discrete version for $\vrh$ satisfying (\ref{n2}) reads:
\begin{equation} \label{r3}
\begin{split}
&\intO{ \Dh (\vr_h b_h) \chi'(b_h) \Phi } - \sum_{ \sigma \in \Sigma_{int} } \intSh{ {Up}[\vr_h b_h, \vu_h] [[\chi' (b_h) \Phi]] }\\
&= \intO{ \Dh \vr_h \chi(b_h) \Phi } - \sum_{ \sigma \in \Sigma_{int} } \intSh{ {Up}[\vr_h \chi (b_h), \vu_h] [[ \Phi]] }\\
&+\sum_{ \sigma \in \Sigma_{int} } \intSh{ \mu_h [[\vr_h]]\ [[ \left( \chi(b_h) - \chi'(b_h) b_h \right) \Phi ]] }\\
& +  \sum_{ \sigma \in \Sigma_{int} } \intSh{
 \Phi^{\rm down} \vr_h^{\rm up}
\left( \widetilde{[[\chi(b_h)]]} - \chi'(b^{\rm down}_h) \widetilde{[[ b_h ]]} \right)  |\Ov{\vu_h} \cdot \vc{n} |
 },
\end{split}
\end{equation}
see \cite[Lemma A.1, Section A.2]{FeKaNo}.

\subsection{Discrete entropy balance equation}

\label{k}

We derive a discrete analogue of the entropy balance equation following step by step its derivation in the continuous setting.

\subsubsection{Discrete kinetic energy equation}

The discrete kinetic energy equation is obtained by taking the scalar product of \eqref{b5} with $\vu_h$, or, at the discrete level,
by taking $\bfPhi = \vu_h \Phi$ in (\ref{n3}):
\[
\begin{split}
\Dh \intOh{ \mh \cdot \vu_h \Phi } &- \sum_{ \sigma \in \Sigma_{int} } \intSh{ {\bf F}_h(\mh,\vuh)
\cdot [[\vuh \Phi ]]  } - \sum_{ \sigma \in \Sigma_{int} } \intSh{ \Ov{p_h} \vc{n} \cdot [[ \vu_h \Phi ]] } \\
&= - h^{\alpha - 1} \sum_{ \sigma \in \Sigma_{int} } \intSh{ [[ \vu_h ]] \cdot [[ \vuh \Phi ]] }.
\end{split}
\]
Next, we use relation (\ref{r3}) for $b_h = \vu_h$,  $\chi(|\vuh|) = \frac{1}{2} |\vuh|^2$ to compute
\[
\begin{split}
\Dh &\intOh{ \mh \cdot \vu_h \Phi } - \sum_{ \sigma \in \Sigma_{int} } \intSh{ {\bf Up}[\vc{m}_h, \vuh] \cdot [[ \vuh \Phi]] }\\
&=
\Dh \intOh{ \vr_h \vu_h \cdot \vu_h \Phi } - \sum_{ \sigma \in \Sigma_{int} } \intSh{  {\bf Up}[\vr_h \vu_h, \vuh ] \cdot
[[ \vuh \Phi]] }\\&=
\Dh \intO{ \frac{1}{2} \vr_h |\vu_h|^2 \Phi } - \sum_{ \sigma \in \Sigma_{int} } \intSh{ {\bf Up} \left[ \frac{1}{2} \vr_h |\vu_h|^2 , \vu_h
\right] [[ \Phi]] }\\
&- \sum_{ \sigma \in \Sigma_{int} } \intSh{ \mu_h [[\vr_h]]\ \left[\left[ \frac{1}{2} |\vu_h|^2  \Phi \right]\right] } +  \frac{1}{2} \sum_{ \sigma \in \Sigma_{int} } \intSh{
 \Phi^{\rm down} \vr_h^{\rm up} |\Ov{\vu_h} \cdot \vc{n} |
[[ \vu_h ]]^2   }.
\end{split}
\]
Consequently, summing up the previous two observations we may infer that
\begin{equation} \label{k1}
\begin{split}
\Dh &\intO{ \frac{1}{2} \vr_h |\vu_h|^2 \Phi } - \sum_{ \sigma \in \Sigma_{int} } \intSh{ {\bf Up} \left[ \frac{1}{2} \vr_h |\vu_h|^2 , \vu_h
\right] [[ \Phi]] }  \\
&= - h^{\alpha - 1} \sum_{ \sigma \in \Sigma_{int} } \intSh{ [[ \vu_h ]] \cdot [[ \vuh \Phi ]] }
+ \sum_{ \sigma \in \Sigma_{int} } \intSh{ \Ov{p_h} \vc{n} \cdot [[ \vu_h \Phi ]] }
- \sum_{ \sigma \in \Sigma_{int} } \intSh{ \mu_h [[ \mh ]] [[\vu_h \Phi ]] } \\
&+\sum_{ \sigma \in \Sigma_{int} } \intSh{ \mu_h [[\vr_h]]\ \left[\left[ \frac{1}{2} |\vu_h|^2  \Phi \right]\right] } -  \frac{1}{2} \sum_{ \sigma \in \Sigma_{int} } \intSh{
 \Phi^{\rm down} {\vr_h}^{\rm up}
 |\Ov{\vu_h} \cdot \vc{n} |
[[ \vu_h ]]^2   }.
\end{split}
\end{equation}
Equation (\ref{k1}) is nothing other than the discrete kinetic energy balance associated to the approximate system (\ref{n2}--\ref{n4}).

\subsubsection{Discrete internal energy equation}

The next step is subtracting (\ref{k1}) from the total energy balance (\ref{n4}):
\[
\begin{split}
\Dh&\intOh{  \vr_h e_h \Phi } - \sum_{ \sigma \in \Sigma_{int} } \intSh{ \Big( {Up}[\vr_h e_h, \vu_h] - \mu_h [[ E_h ]]  \Big)
[[\Phi]]  }  \\
&= - h^{\alpha - 1} \sum_{ \sigma \in \Sigma_{int} } \intSh{ [[ \vu_h ]] \cdot [[ \Phi \vuh ]] } +
h^{\alpha - 1} \sum_{ \sigma \in \Sigma_{int} } \intSh{ [[ \vu_h ]]^2 \Ov{\Phi} }\\
& + h^{\alpha - 1} \sum_{ \sigma \in \Sigma_{int} } \intSh{ [[ \vu_h ]] \cdot [[ \vuh \Phi ]] }
- \sum_{ \sigma \in \Sigma_{int} } \intSh{ \Ov{p_h \Phi} \vc{n} \cdot [[ \vu_h ]] }
+ \sum_{ \sigma \in \Sigma_{int} } \intSh{ \mu_h [[ \mh ]] [[\vu_h \Phi ]] } \\
&-\sum_{ \sigma \in \Sigma_{int} } \intSh{ \mu_h [[\vr_h]]\ \left[\left[ \frac{1}{2} |\vu_h|^2  \Phi \right]\right] } +  \frac{1}{2} \sum_{ \sigma \in \Sigma_{int} } \intSh{  \Phi^{\rm down} {\vr_h}^{\rm up} |\Ov{\vu_h} \cdot \vc{n} |
[[ \vu_h ]]^2   },
\end{split}
\]
or, reordered,
\[
\begin{split}
\Dh &\intOh{ \vr_h e_h \Phi } - \sum_{ \sigma \in \Sigma_{int} } \intSh{ \Big( {Up}[\vr_h e_h, \vu_h] - \mu_h [[ \vr_h e_h ]]  \Big)
[[\Phi]]  }  \\
&=  h^{\alpha - 1} \sum_{ \sigma \in \Sigma_{int} } \intSh{ [[ \vu_h ]]^2 \Ov{\Phi} }
+ \frac{1}{2} \sum_{ \sigma \in \Sigma_{int} } \intSh{  \Phi^{\rm down} {\vr_h}^{\rm up} |\Ov{\vu_h} \cdot \vc{n} |
[[ \vu_h ]]^2   }
- \sum_{ \sigma \in \Sigma_{int} } \intSh{ \Ov{p_h \Phi} [[ \vu_h ]] \cdot \vc{n}  }\\
&+ \sum_{ \sigma \in \Sigma_{int} } \intSh{ \mu_h [[ \vr_h \vu_h ]] [[\vu_h \Phi ]] } -\sum_{ \sigma \in \Sigma_{int} } \intSh{ \mu_h [[\vr_h]]\ \left[\left[ \frac{1}{2} |\vu_h|^2  \Phi \right]\right] }  \\
& - \sum_{ \sigma \in \Sigma_{int} } \intSh{ \mu_h \left[\left[ \frac{1}{2} \vr_h |\vu_h|^2 \right]\right] [[\Phi]] }.
\end{split}
\]
Finally, using the product rule, we obtain
\[
\begin{split}
[[ \vr_h \vu_h ]] [[\vu_h \Phi ]] &- \frac{1}{2} [[\vr_h]]\ [[ |\vu_h|^2  \Phi ]] -
\frac{1}{2} [[ \vr_h |\vu_h|^2 ]] [[\Phi]]\\ 
&= \Ov{\vr_h} [[ \vu_h ]] \cdot [[\vu_h ]] \Ov{\Phi}  + \Ov{\vr_h}\ \Ov{\vu_h} \cdot [[ \vu_h ]] [[\Phi]] \\
&\hskip 3.38cm +
\frac{1}{2} [[ \vr_h ]] \Ov{\vu}_h  \cdot [[\vu_h \Phi ]]
- \frac{1}{2} [[\vr_h]] \Ov{\vu_h} \cdot [[\vu_h]] \Ov{\Phi}  -
\frac{1}{2} [[ \vr_h |\vu_h|^2 ]] [[\Phi]] \\
&= \Ov{\vr_h} [[ \vu_h ]] \cdot [[\vu_h ]] \Ov{\Phi}  +
\Ov{\vr_h}\ \Ov{\vu_h} \cdot [[ \vu_h ]] [[\Phi]] + \frac{1}{2} [[ \vr_h]] |\Ov{\vu_h}|^2 [[\Phi]]-
\frac{1}{2} [[ \vr_h |\vu_h|^2 ]] [[\Phi]] \\
&=\Ov{\vr_h} [[ \vu_h ]] \cdot [[\vu_h ]] \Ov{\Phi} +\Ov{\vr_h}\ \Ov{\vu_h} \cdot [[ \vu_h ]] [[\Phi]]
- \frac{1}{2} \Ov{\vr_h} [[ \vu_h \cdot \vu_h ]] [[ \Phi ]] = \Ov{\vr_h} [[\vu_h ]]^2 \Ov{\Phi}.
\end{split}
\]
Consequently, we record the internal energy balance in the form
\begin{equation} \label{k2}
\begin{aligned}
 \Dh\intOh{  \vr_h e_h \Phi } &- \sum_{ \sigma \in \Sigma_{int} } \intSh{ \Big( {Up}[\vr_h e_h, \vu_h] - \mu_h [[ \vr_h e_h ]]  \Big)
[[\Phi]]  }  \\
&=  h^{\alpha - 1} \sum_{ \sigma \in \Sigma_{int} } \intSh{ [[ \vu_h ]]^2 \Ov{\Phi} }
+ \frac{1}{2} \sum_{ \sigma \in \Sigma_{int}} \intSh{  \Phi^{\rm down} {\vr_h}^{\rm up}|\Ov{\vu_h} \cdot \vc{n} |
[[ \vu_h ]]^2   }  \\
&+ \sum_{ \sigma \in \Sigma_{int} } \intSh{ \mu_h \Ov{\vr_h} [[\vu_h ]]^2 \Ov{\Phi} }- \sum_{ \sigma \in \Sigma_{int} } \intSh{ \Ov{p_h \Phi} [[ \vu_h ]] \cdot \vc{n}  }.
\end{aligned}
\end{equation}

\subsubsection{Discrete entropy balance}

At this stage, we are ready to derive the discrete entropy balance together with its renormalization.
Dividing equation (\ref{k2}) on $\vt_h$, we get
\begin{align*}
c_v&\intOh{ \Dh (\vr_h \vt_h) \left(\frac{\Phi}{\vt_h} \right) } - c_v \sum_{ \sigma \in \Sigma_{int} } \intSh{ {Up}[\vr_h \vt_h, \vu_h]
\left[\left[\frac{\Phi}{\vt_h} \right]\right]  }  \\
&=  h^{\alpha - 1} \sum_{ \sigma \in \Sigma_{int} } \intSh{ [[ \vu_h ]]^2 \Ov{ \left( \frac{\Phi}{\vt_h} \right) }  }
+ \frac{1}{2} \sum_{ \sigma \in \Sigma_{int} } \intSh{  \left( \frac{\Phi}{\vt_h} \right)^{\rm down}  \vr_h^{\rm up}
|\Ov{\vu_h} \cdot \vc{n} |
[[ \vu_h ]]^2   }\\
&+ \sum_{ \sigma \in \Sigma_{int} } \intSh{ \mu_h \Ov{\vr_h} [[\vu_h ]]^2 \Ov{ \left( \frac{\Phi}{\vt_h} \right) } }
- \sum_{ \sigma \in \Sigma_{int} } \intSh{ [[ \vu_h ]] \cdot \vc{n} \Ov{ \vr_h \Phi }  }
- c_v \sum_{ \sigma \in \Sigma_{int} } \intSh{ \mu_h [[ \vr_h \vt_h ]] \ \left[\left[ \frac{\Phi}{\vt_h} \right]\right] }.
\end{align*}
Next, by virtue of formula (\ref{r3}),
\[
\begin{split}
c_v& \intOh{ \Dh (\vr_h \vt_h) \left(\frac{\Phi}{\vt_h} \right) } - c_v \sum_{ \sigma \in \Sigma_{int} } \intSh{ {Up}[\vr_h \vt_h, \vu_h]
\left[\left[\frac{\Phi}{\vt_h} \right]\right]  }\\
&= \Dh \intO{ \vr_h \log(\vt^{c_v}_h) \Phi } - \sum_{ \sigma \in \Sigma_{int} } \intSh{ {\bf Up}[\vr_h \log(\vt_h^{c_v}) , \vu_h] [[ \Phi]] }\\
&+ c_v \sum_{ \sigma \in \Sigma_{int} } \intSh{ \mu [[\vr_h]]\ [[ \left( \log(\vt_h) - 1 \right) \Phi ]] }\\
& +  c_v \sum_{ \sigma \in \Sigma_{int} } \intSh{  \Phi^{\rm down} {\vr_h}^{\rm up}
\left( \widetilde{[[\log(\vt_h)]]} - \frac{1}{{\vt_h^{\rm down}}} \widetilde{[[ \vt_h ]]} \right) |\Ov{\vu_h} \cdot \vc{n} |}  .
\end{split}
\]
Consequently,
\begin{equation} \label{k3}
\begin{aligned}
\Dh &\intOh{ \vr_h \log(\vt^{c_v}_h) \Phi } - \sum_{ \sigma \in \Sigma_{int} } \intSh{ {Up}[\vr_h \log(\vt_h^{c_v}) , \vu_h] [[ \Phi]] } \\
&=  h^{\alpha - 1} \sum_{ \sigma \in \Sigma_{int} } \intSh{ [[ \vu_h ]]^2 \Ov{ \left( \frac{\Phi}{\vt_h} \right) }  }
+ \frac{1}{2} \sum_{ \sigma \in \Sigma_{int} } \intSh{ \left( \frac{\Phi}{\vt_h} \right)^{\rm down}  \vr_h^{\rm up}  |\Ov{\vu_h} \cdot \vc{n} |
[[ \vu_h ]]^2   }\\
&+ \sum_{ \sigma \in \Sigma_{int} } \intSh{ \mu_h \Ov{\vr_h} [[\vu_h ]]^2 \Ov{ \left( \frac{\Phi}{\vt_h} \right) } }
- \sum_{ \sigma \in \Sigma_{int} } \intSh{ [[ \vu_h ]] \cdot \vc{n} \Ov{ \vr_h \Phi }  }
\\ &- c_v \sum_{ \sigma \in \Sigma_{int} } \intSh{ \mu_h [[ \vr_h \vt_h ]] \ \left[\left[ \frac{\Phi}{\vt_h} \right]\right] }
- c_v \sum_{ \sigma \in \Sigma_{int} } \intSh{ \mu_h [[\vr_h]]\ [[ \left( \log(\vt_h) - 1 \right) \Phi ]] }\\
& -  c_v \sum_{ \sigma \in \Sigma_{int} } \intSh{   \Phi^{\rm down} {\vr_h}^{\rm up}
\left( \widetilde{[[\log(\vt_h)]]} - \frac{1}{{\vt_h^{\rm down}}} \widetilde{[[ \vt_h ]]} \right) |\Ov{\vu_h} \cdot \vc{n} |}.
\end{aligned}
\end{equation}
Finally, we
consider $b(\vr) = \vr \log(\vr)$ in the renormalized equation \eqref{r1}:
\begin{equation} \label{k4}
\begin{split}
\Dh &\intOh{ \vr_h \log(\vr_h) \Phi } - \sum_{ \sigma \in \Sigma_{int} } \intSh{  {Up}[\vr_h \log (\vr_h), \vu_h]
 [[\Phi]] }
\\ &= - \sum_{ \sigma \in \Sigma_{int} } \intSh{ \mu_h [[\vrh]] [[ b'(\vrh) \Phi ]] }
- \sum_{ \sigma \in \Sigma_{int} } \intSh{
 \Phi^{\rm down}
\Big( \widetilde{[[ b(\vrh) ]]} - b'( \vrh^{\rm down} ) \widetilde{[[ \vrh ]]} \Big)  | \Ov{\vuh} \cdot \vc{n} | }
\\
&- \sum_{ \sigma \in \Sigma_{int} } \intSh{ [[{\vuh} ]] \cdot \vc{n}  \Ov{ \vr_h \Phi } }.
\end{split}
\end{equation}
Subtracting (\ref{k4}) from (\ref{k3}) and introducing the entropy
$ \displaystyle
s_h = \log \left( \frac{\vt^{c_v}_h}{\vrh} \right)
$
we obtain
\begin{align}
 \Dh&\intOh{ \vr_h s_h \Phi } - \sum_{ \sigma \in \Sigma_{int} } \intSh{ {Up}[\vr_h s_h , \vu_h] [[ \Phi]] } \nonumber  \\
&=  h^{\alpha - 1} \sum_{ \sigma \in \Sigma_{int} } \intSh{ [[ \vu_h ]]^2 \Ov{ \left( \frac{\Phi}{\vt_h} \right) }  }
+ \frac{1}{2} \sum_{ \sigma \in \Sigma_{int} } \intSh{ \left( \frac{\Phi}{\vt_h} \right)^{\rm down}  \vr_h^{\rm up} |\Ov{\vu_h} \cdot \vc{n} |
[[ \vu_h ]]^2   }\nonumber \\
&+ \sum_{ \sigma \in \Sigma_{int} } \intSh{ \mu_h \Ov{\vr_h} [[\vu_h ]]^2 \Ov{ \left( \frac{\Phi}{\vt_h} \right) } }
+ \sum_{ \sigma \in \Sigma_{int} } \intSh{  \Phi^{\rm down}
\Big( \widetilde{[[ b(\vrh) ]]} - b'( \vrh^{\rm down} ) \widetilde{[[ \vrh ]]} \Big)  | \Ov{\vuh} \cdot \vc{n} |  }
 \nonumber \\
&-  c_v \sum_{ \sigma \in \Sigma_{int} } \intSh{  \Phi^{\rm down} {\vr_h}^{\rm up}
\left( \widetilde{[[\log(\vt_h)]]} - \frac{1}{{\vt_h^{\rm down}}} \widetilde{[[ \vt_h ]]} \right) |\Ov{\vu_h} \cdot \vc{n} |}
 \label{k5} \\ &- c_v \sum_{ \sigma \in \Sigma_{int} } \intSh{ \mu_h [[ \vr_h \vt_h ]] \ \left[\left[ \frac{\Phi}{\vt_h} \right]\right] }
- c_v \sum_{ \sigma \in \Sigma_{int} } \intSh{ \mu_h [[\vr_h]]\ [[ \left( \log(\vt_h) - 1 \right) \Phi ]] }
\nonumber \\
&+\sum_{ \sigma \in \Sigma_{int} } \intSh{ \mu_h [[\vrh]] [[ b'(\vrh) \Phi ]] } ,\
\mbox{ where }\ b(\vr) = \vr \log(\vr).\nonumber
\end{align}
This is the physical entropy balance associated to (\ref{n2}--\ref{n4}). At this stage, it is not obvious how to handle the last three integrals in (\ref{k5}), however, this will be fixed in the following section.

\subsubsection{Entropy renormalization}

Consider $\chi$ - a non--decreasing, concave, twice continuously differentiable function on $R$ that is bounded from above.
Applying formula \eqref{r3} in \eqref{k5} we get
\begin{align*}
\Dh&\intOh{  \vr_h \chi(s_h) \Phi } - \sum_{ \sigma \in \Sigma_{int} } \intSh{ {Up}[\vr_h \chi(s_h) , \vu_h] [[ \Phi]] } \\
&=  h^{\alpha - 1} \sum_{ \sigma \in \Sigma_{int} } \intSh{ [[ \vu_h ]]^2 \Ov{ \left( \frac{\chi'(s_h)\Phi}{\vt_h} \right) }  }
+ \frac{1}{2} \sum_{ \sigma \in \Sigma_{int} } \intSh{   \left( \frac{\chi'(s_h) \Phi}{\vt_h} \right)^{\rm down} \vr_h^{\rm up} |\Ov{\vu_h} \cdot \vc{n} |
[[ \vu_h ]]^2 }  \\
&+ \sum_{ \sigma \in \Sigma_{int} } \intSh{ \mu_h \Ov{\vr_h} [[\vu_h ]]^2 \Ov{ \left( \frac{\chi'(s_h) \Phi}{\vt_h} \right) } } \\
&
+ \sum_{ \sigma \in \Sigma_{int} } \intSh{ \left(\chi'(s_h) \Phi\right)^{\rm down}  \left( \widetilde{[[ b(\vrh) ]]} -
b'( \vr_h^{\rm down}  ) \widetilde{[[ \vrh ]]} \right)| \Ov{\vuh} \cdot \vc{n} |  } \\
&-  c_v \sum_{ \sigma \in \Sigma_{int} } \intSh{  \left(\chi'(s_h)\Phi\right)^{\rm down} \vr_h^{\rm up}
\left( \widetilde{[[\log(\vt_h)]]} - \frac{1}{\vt_h^{\rm down}} \widetilde{[[ \vt_h ]]} \right)|\Ov{\vu_h} \cdot \vc{n} |   }\\
&-  \sum_{ \sigma \in \Sigma_{int} } \intSh{  \Phi^{\rm down} \vr_h^{\rm up}
\left( \widetilde{[[\chi(s_h)]]} - \chi'(s^{\rm down}_h) \widetilde{[[ s_h ]]} \right)  |\Ov{\vu_h} \cdot \vc{n} | }
\\ &- c_v \sum_{ \sigma \in \Sigma_{int} } \intSh{ \mu_h [[ \vr_h \vt_h ]]  \left[\left[ \frac{\chi'(s_h)\Phi}{\vt_h} \right]\right] }
- c_v \sum_{ \sigma \in \Sigma_{int} } \intSh{ \mu_h [[\vr_h]]\ [[ \left( \log(\vt_h) - 1 \right) \chi'(s_h) \Phi ]] }
\\
&+\sum_{ \sigma \in \Sigma_{int} } \intSh{ \mu_h [[\vrh]]\ [[ b'(\vrh) \chi'(s_h)\Phi ]] } \\
&- \sum_{ \sigma \in \Sigma_{int} } \intSh{  \mu_h [[\vr_h]]\ [[ \left( \chi(s_h) - \chi'(s_h) s_h \right) \Phi ]] },
\ \mbox{ where }\ b(\vr) = \vr \log(\vr).
\end{align*}
Next, we compute
\begin{align*}
&- c_v \sum_{ \sigma \in \Sigma_{int} } \intSh{ \mu_h [[ \vr_h \vt_h ]] \left[\left[ \frac{\chi'(s_h)\Phi}{\vt_h} \right]\right] }
- c_v \sum_{ \sigma \in \Sigma_{int} } \intSh{ \mu_h [[\vr_h]]\ [[ \left( \log(\vt_h) - 1 \right) \chi'(s_h) \Phi ]] }
\\
&+\sum_{ \sigma \in \Sigma_{int} } \intSh{ \mu_h [[\vrh]]\ [[ b'(\vrh) \chi'(s_h)\Phi ]] } - \sum_{ \sigma \in \Sigma_{int} } \intSh{ \mu_h [[\vr_h]]\ \left[\left[ \left( \chi(s_h) - \chi'(s_h) s_h \right) \Phi \right]\right] }\\
&=- c_v \sum_{ \sigma \in \Sigma_{int} } \intSh{ \mu_h [[ \vr_h \vt_h ]] \left[\left[ \frac{\chi'(s_h)\Phi}{\vt_h} \right]\right] }
- c_v \sum_{ \sigma \in \Sigma_{int} } \intSh{ \mu_h [[\vr_h]]\ [[ \log(\vt_h) \chi'(s_h) \Phi ]] }
\\
&\hskip 0.8cm +\sum_{ \sigma \in \Sigma_{int} } \intSh{ \mu_h [[\vrh]] \ [[ \log(\vrh) \chi'(s_h)\Phi ]] } - \sum_{ \sigma \in \Sigma_{int} } \intSh{ \mu_h [[\vr_h]]\ [[ \left( \chi(s_h) - \chi'(s_h) s_h \right) \Phi ]] }\\
&\hskip 0.8cm  +(c_v + 1) \sum_{ \sigma \in \Sigma_{int} } \intSh{ \mu_h [[ \vrh ]] \ [[\chi'(s_h) \Phi ]] }\\
&=- c_v \sum_{ \sigma \in \Sigma_{int} } \intSh{ \mu_h [[ \vr_h \vt_h ]] \left[\left[ \vrh \frac{\chi'(s_h)\Phi}{\vrh \vt_h} \right]\right] }
- \sum_{ \sigma \in \Sigma_{int} } \intSh{  \mu_h [[\vr_h]] \left[\left[ \Big( \chi(s_h)  - (c_v + 1)\chi'(s_h) \Big) \Phi \right]\right] }\\
&=- c_v \sum_{ \sigma \in \Sigma_{int} } \intSh{ \mu_h [[ p_h ]]  \left[\left[ \vrh \frac{\chi'(s_h)\Phi}{p_h} \right]\right] }
- \sum_{ \sigma \in \Sigma_{int} } \intSh{ \mu_h [[\vr_h]] \left[\left[ \Big( \chi(s_h) - (c_v + 1)\chi'(s_h) \Big) \Phi \right]\right] }\\ &=
- \sum_{ \sigma \in \Sigma_{int} } \intSh{ \mu_h [[ \Phi \nabla_{\vr} (\vr_h \chi(s_h)) ]] \ [[\vr_h ]] }
- \sum_{ \sigma \in \Sigma_{int} } \intSh{ \mu_h [[ \Phi \nabla_p (\vr_h \chi(s_h)) ]] \ [[ p_h ]] }.
\end{align*}
Thus we infer with the general entropy inequality
\begin{equation} \label{k6}
\begin{aligned}
 \Dh&\intOh{  \vr_h \chi(s_h) \Phi } - \sum_{ \sigma \in \Sigma_{int} } \intSh{  Up[\vr_h \chi(s_h) , \vu_h] [[ \Phi]] } \\
&=  h^{\alpha - 1} \sum_{ \sigma \in \Sigma_{int} } \intSh{ [[ \vu_h ]]^2 \Ov{ \left( \frac{\chi'(s_h)\Phi}{\vt_h} \right) }  }
+ \frac{1}{2} \sum_{ \sigma \in \Sigma_{int} } \intSh{  \left( \frac{\chi'(s_h) \Phi}{\vt_h} \right)^{\rm down} \vr_h^{\rm up}  |\Ov{\vu_h} \cdot \vc{n} |
[[ \vu_h ]]^2   }\\
&+ \sum_{ \sigma \in \Sigma_{int} } \intSh{ \mu_h \Ov{\vr_h} [[\vu_h ]]^2 \Ov{ \left( \frac{\chi'(s_h) \Phi}{\vt_h} \right) } } \\
&+ \sum_{ \sigma \in \Sigma_{int} } \intSh{ \left(\chi'(s_h) \Phi\right)^{\rm down}  \left( \widetilde{[[ b(\vrh) ]]} -
b'( \vr_h^{\rm down}  ) \widetilde{[[ \vrh ]]} \right)| \Ov{\vuh} \cdot \vc{n} | } \\
&-  c_v \sum_{ \sigma \in \Sigma_{int} } \intSh{ \left(\chi'(s_h)\Phi\right)^{\rm down} \vr_h^{\rm up}
\left( \widetilde{[[\log(\vt_h)]]} - \frac{1}{\vt_h^{\rm down}} \widetilde{[[ \vt_h ]]} \right)|\Ov{\vu_h} \cdot \vc{n} |    }\\
&-  \sum_{ \sigma \in \Sigma_{int} } \intSh{  \Phi^{\rm down} \vr_h^{\rm up}
\left( \widetilde{[[\chi(s_h)]]} - \chi'(s^{\rm down}_h) \widetilde{[[ s_h ]]} \right) |\Ov{\vu_h} \cdot \vc{n} |  } \\
&- \sum_{ \sigma \in \Sigma_{int} } \intSh{ \mu_h [[ \Phi \nabla_{\vr} (\vr_h \chi(s_h)) ]] \ [[\vr_h ]] }
- \sum_{ \sigma \in \Sigma_{int} } \intSh{ \mu_h [[ \Phi \nabla_p (\vr_h \chi(s_h)) ]]\ [[ p_h ]] },\
b(\vr) = \vr \log(\vr).
\end{aligned}
\end{equation}
 Note that the last two integrals in (\ref{k6}) can be rewritten using the product rule as
\begin{equation}
\begin{aligned}\label{en_num_flux_muh}
&- \sum_{ \sigma \in \Sigma_{int} } \intSh{ \mu_h [[\Phi]]\left( \Ov{\nabla_{\vr} (\vr_h \chi(s_h))}[[\vr_h ]]+\Ov{\nabla_p (\vr_h \chi(s_h))}[[p_h ]]\right)  } \\
&
- \sum_{ \sigma \in \Sigma_{int} } \intSh{ \mu_h \Ov{\Phi}[[\nabla_{\vr} (\vr_h \chi(s_h))]]\ [[\vrh]]+ [[  \nabla_p (\vr_h \chi(s_h)) ]]\ [[ p_h ]] }.
\end{aligned}
\end{equation}
The first sum in \eqref{en_num_flux_muh}  together with the upwind term in \eqref{k6},
\begin{equation}
\begin{aligned}\label{num_flux_entropy}
- \sum_{ \sigma \in \Sigma_{int} } \intSh{ \left( Up[\vr_h \chi(s_h) , \vu_h]   +\mu_h \left( \Ov{\nabla_{\vr} (\vr_h \chi(s_h))}[[\vr_h ]]+\Ov{\nabla_p (\vr_h \chi(s_h))}[[p_h ]]\right)  \right)[[ \Phi]]},
\end{aligned}
\end{equation}
represent the numerical  entropy flux.
The rest in \eqref{k6} and \eqref{en_num_flux_muh} gives the numerical entropy production, cf. \cite{FLM18,FjMiTa1,FjMiTa2}.
Recall that the total entropy
\[
(\vr,p) \mapsto -\vr \chi (s(\vr,p)) = -\vr \chi \left( \log\left( \frac{\vt^{c_v}}{\vr} \right) \right) =-
\vr \chi \left( \frac{1}{\gamma - 1} \log \left( \frac{p}{\vr^\gamma} \right) \right)
\]
 is a convex function of the variables $\vr$ and $p.$ In particular, $-\nabla_{\vr,p} (\vr\chi(s(\vr_h,p_h)))$ is monotone, and therefore the term in the second line of \eqref{en_num_flux_muh} is non--negative.
It is worthwhile to mention that the discrete entropy inequality \eqref{k6} is a discrete version of \eqref{entropy_cont} with $\kappa=c_vh\vr\lambda,$ $\lambda=\frac{1}{2}|\Ov{\vu_h}\cdot\vn|+\mu_h.$

\section{Stability}
\label{AS}

Having established all necessary ingredients, we are ready to discuss the available {\it a priori} bounds for solutions of the semi--discrete scheme (\ref{n2}--\ref{n4}).

\subsection{Mass and energy conservation}

Taking $\Phi \equiv 1$ in the equation of continuity (\ref{n2}) yields the total mass conservation
\begin{equation} \label{as1}
\intOh{ \vr_h (t, \cdot) } = \intOh{ \vr_{0,h} } = M_0 > 0,\ t \geq 0.
\end{equation}
A similar argument applied to the total energy balance yields
\begin{equation} \label{as2}
\intOh{ E_h (t, \cdot) } = \intOh{E_{0,h}} = E_0 > 0,\ t \geq 0.
\end{equation}

\subsection{Minimum entropy principle}

An important source of {\it a priori} bounds is the minimum entropy principle that can be derived from the entropy balance with the choice
\[
\Phi = 1,\ \chi(s) = |s - s_0|^{-},\
- \infty < s_0 < \min s_h(0).
\]
As
\[
\vr \mapsto \vr \log(\vr) \ \mbox{is convex},\
\vt \mapsto \log(\vt) \ \mbox{concave},\ s \mapsto \chi(s) \ \mbox{concave},\
\mbox{and}\ (\vr, p) \mapsto \vr \chi (s(\vr, p)) \ \mbox{concave,}
\]
all integrals on the right--hand side of (\ref{k6}) are non--negative, and we may infer that
\[
\intOh{ \vr_h(t) | s_h(t) - s_0 |^- } \geq 0 \ \mbox{ for any }\ t \geq 0.
\]
Consequently, we have obtained the minimum entropy principle
\begin{equation} \label{as3}
s_h(t) \geq s_0 \ \mbox{ for all }\ t \geq 0.
\end{equation}

\subsection{Positivity of the pressure, existence of the temperature}

The entropy as a function of $\vr$ and $p$ reads
\[
s = \frac{1}{\gamma -1 }\log \left( \frac{p}{\vr^\gamma} \right);
\]
whence it follows immediately from (\ref{as3}) that
\begin{equation} \label{as4}
0 < \exp \{(\gamma - 1) s_0 \} \leq \frac{p_h(t) }{\vr^\gamma_h(t)} \ \mbox{ for all }\ t \geq 0.
\end{equation}
In particular, the pressure is positive as long as the density is positive, and we may set
\[
\vt_h(t) = \frac{p_h(t)}{\vr_h(t)}.
\]
Evoking the energy bound (\ref{as2}) we get
\begin{equation} \label{as5}
\frac{1}{2} \intOh{ \frac{|\vc{m}_h(t)|^2}{\vrh(t)} } +
c_v \intOh{ \vr_h(t) \vt_h(t) } \leq E_0 \ \mbox{ for all } \ t \geq 0.
\end{equation}
Thus going back to (\ref{as4}) we obtain
\begin{equation} \label{as6}
\intOh{ \vr_h^\gamma(t) } \aleq \intOh{ p_h(t) } \aleq E_0 \ \mbox{ for all }\ t \geq 0.
\end{equation}

\subsection{Positivity of the density}

The crucial property for the approximate scheme to be valid is positivity of the density $\vrh$ at least at the discrete level, meaning for any $h > 0$.  We will show that, for any $T>0,$ there exists $\underline{\vr}=\underline{\vr}(h,T)>0,$ such that $\vrh(t)\geq \underline{\vr}>0$ for all $t\in [0,T].$  To see this, we first evoke the kinetic energy balance (\ref{k1}) with $\Phi = 1$. Seeing that
\[
- \sum_{ \sigma \in \Sigma_{int} } \intSh{ \mu_h [[ \mh ]] [[\vu_h ]] }
+\sum_{ \sigma \in \Sigma_{int} } \intSh{ \mu_h [[\vr_h]]\ \left[\left[ \frac{1}{2} |\vu_h|^2   \right]\right] }
= - \sum_{ \sigma \in \Sigma_{int} } \intSh{ \mu_h \Ov{\vrh} [[\vu_h ]]^2 },
\]
we may integrate (\ref{k1}) in time and use the energy bound (\ref{as5}) to deduce
\[
\begin{split}
& h^{\alpha - 1} \int_0^T \sum_{ \sigma \in \Sigma_{int} } \intSh{ [[ \vu_h ]]^2  } \dt +
\int_0^T \sum_{ \sigma \in \Sigma_{int} } \intSh{ \mu_h \Ov{\vrh} [[\vu_h ]]^2 } \dt \\
&+\frac{1}{2} \int_0^T \sum_{ \sigma \in \Sigma_{int} } \intSh{  \vr_h^{\rm up} |\Ov{\vu_h} \cdot \vc{n} |
[[ \vu_h ]]^2   } \dt
\aleq \left(1 + \sum_{ \sigma \in \Sigma_{int} } \int_0^T \intSh{ \Ov{p_h} \vc{n} \cdot [[ \vu_h ]] } \right) \dt.
\end{split}
\]
Finally, we again use (\ref{as5}) combined with the negative $L^p$--estimates (\ref{n4b}) and H\" older's inequality to conclude
\begin{equation} \label{as7}
\int_0^T \sum_{ \sigma \in \Sigma_{int} } \intSh{ [[ \vu_h ]]^2  } \dt \aleq \omega(h),
\end{equation}
where $\omega(h)$ denotes a generic function that may blow up in the asymptotic regime $h \to 0$. In particular, relation (\ref{as7}) implies
\begin{equation} \label{as8}
\int_0^T \left( \sup_{ \sigma \in \Sigma_{int} } [[ \vu_h ]]^2 \right) \dt \aleq \omega(h),
\end{equation}
with another $\omega(h)$ generally different from its counterpart in (\ref{as7}).

Next, we revisit the renormalized equation of continuity (\ref{r1}), again with $\Phi = 1$, obtaining
\[
\intOh{ \Dt b(\vr_h) }
+ \sum_{ \sigma \in \Sigma_{int} } \intSh{ [[ \Ov{\vuh} ]] \cdot \vc{n} \Ov{ \Big( b(\vr_h) - b'(\vr_h) \vr_h \Big) } }
\leq 0
\]
for any convex $b$. Thus the specific choice $b(\vr) = |\vr - \underline{\vr} |^-$ gives rise to the inequality
\[
\Dt \intOh{ |\vr_h - \underline{\vr} |^- }
+ \underline{\vr} \sum_{ \sigma \in \Sigma_{int} } \intSh{ [[ \Ov{\vuh} ]] \cdot \vc{n} \ \Ov{ 1_{\vrh(t) \leq \underline{\vr}} } }
\leq 0.
\]
In view of (\ref{as8}), we can find a positive constant $\underline{\vr} = \underline{\vr}(h,T) > 0$ small enough so that
\[
\intOh{ |\vr_h(t) - \underline{\vr} |^- } < 0 \ \mbox{ for all }\ t \in [0,T].
\]
In other words
\begin{equation} \label{as9}
\vr_h (t) \geq \underline{\vr}(h,T) > 0 \ \mbox{ for all }\ t \in [0,T].
\end{equation}

\begin{Remark} \label{Rn1}

Of course, the estimate (\ref{as9}) is not uniform, neither with respect to $T$ nor for $h \to 0$. In particular, the
asymptotic limit may experience vacuum zone where the density vanishes.

\end{Remark}

\subsection{Existence of approximate solutions}

Having established positivity of the density on any compact time interval, we have closed the {\it a priori} bounds that
guarantee global existence for the semi--discrete system at any level $h > 0$.

\begin{Theorem} \label{Tm1}

Suppose that the initial data $\vr_{0,h}$, $\vm_{0,h}$, $E_{0,h}$ satisfy
\[
\vr_{0,h} \geq \underline{\vr} > 0,\ E_{0,h} - \frac{1}{2} \frac{ |\vc{m}_{0,h}|^2 }{\vr_{0,h}} > 0.
\]
\\
Then the semi--discrete approximate system (\ref{n2}--\ref{n4}) admits a unique global-in-time solution
$[\vrh, \vm_h, E_h]$ such that
\[
\vrh(t) > 0, \ E_{h}(t) - \frac{1}{2} \frac{ |\vm_h(t) |^2 }{\vrh(t)} > 0 \ \mbox{ for any }\ t \geq 0.
\]
\\
Moreover, the renormalized entropy balance (\ref{k6}) holds.

\end{Theorem}

\subsection{Entropy estimates}

We close this section by showing the uniform bounds provided by the dissipation mechanism hidden in the entropy production rate.
We start by observing that
\begin{equation} \label{as10}
\intOh{ \vr_h s_h (t) } \aleq 1 + \intOh{ E_h(t) } \leq 1 + E_0.
\end{equation}
Indeed, in view of the minimum entropy principle established in (\ref{as4}), it is enough to observe that
\[
\vr_h \log \left( \frac{{\vt}_h^{c_v} }{\vrh} \right) \aleq 1 + \vrh \vt_h \ \mbox{ provided }\ 0 < \vrh \aleq \vt_h^{c_v}.
\]
Seeing that $\vrh \log (\vrh)$ is controlled by (\ref{as6}) we restrict ourselves to
$\vrh \log (\vt^{c_v}_h )$. Here,
\[
\vrh \log (\vt^{c_v}_h ) \aleq \vrh \vt_h \aleq E_0 \ \mbox{ if }\ \vt_h \geq 1,
\]
while
\[
|\vrh \log (\vt^{c_v}_h )| \leq \vt^{c_v}_h |\log (\vt^{c_v}_h) | \aleq 1
\ \mbox{ for }\ \vt_h \leq 1.
\]
Thus we have shown (\ref{as10}).

In accordance with (\ref{as10}), we can take
 $\Phi = 1$, $\chi_{\varepsilon}(s) = \min\{s,\frac{1}{\varepsilon}\}$ in the renormalized entropy balance \eqref{k6}. Letting $\varepsilon \rightarrow 0$  we obtain the uniform estimate:
\begin{equation} \label{as11}
\begin{split}
h^{\alpha - 1} &\int_0^T \sum_{ \sigma \in \Sigma_{int} } \intSh{ [[ \vu_h ]]^2 \Ov{ \left( \frac{1}{\vt_h} \right) }  } \dt
+ \frac{1}{2}  \int_0^T \sum_{ \sigma \in \Sigma_{int} } \intSh{   \left( \frac{1}{\vt_h} \right)^{\rm down} \vr_h^{\rm up}
 |\Ov{\vu_h} \cdot \vc{n} |
[[ \vu_h ]]^2   } \dt \\
+& \int_0^T \sum_{ \sigma \in \Sigma_{int} } \intSh{ \mu_h \Ov{\vr_h} [[\vu_h ]]^2 \Ov{ \left( \frac{1}{\vt_h} \right) } } \dt
+ \int_0^T \sum_{ \sigma \in \Sigma_{int} } \intSh{
 \left( \widetilde{[[ b(\vr) ]]} -
b'( {\vr_h^{\rm down} } ) \widetilde{[[ \vrh ]]} \right)| \Ov{\vuh} \cdot \vc{n} |  } \dt
 \\
-&  c_v \int_0^T \sum_{ \sigma \in \Sigma_{int} } \intSh{
 \vr_h^{\rm down}
\left( \widetilde{[[\log(\vt_h)]]} - \frac{1}{\vt_h^{\rm down}} \widetilde{ [[ \vt_h ]]} \right)  |\Ov{\vu_h} \cdot \vc{n} |  } \dt \\
-& \int_0^T \sum_{ \sigma \in \Sigma_{int} } \intSh{ \mu_h [[ \nabla_{\vr} (\vr_h s_h) ]]\  [[\vr_h ]] } \dt
- \int_0^T \sum_{ \sigma \in \Sigma_{int} } \intSh{ \mu_h [[ \nabla_p (\vr_h s_h) ]] \ [[ p_h ]] } \dt \aleq (1 + E_0),
\end{split}
\end{equation}
where $b(\vr) = \vr \log(\vr)$.
As for the last two integrals in (\ref{as11}), we can check by direct manipulation that
\[
\begin{split}
- \int_0^T &\sum_{ \sigma \in \Sigma_{int} } \intSh{ \mu_h [[ \nabla_{\vr} (\vr_h s_h) ]]\ [[\vr_h ]] } \dt
- \int_0^T \sum_{ \sigma \in \Sigma_{int} } \intSh{ \mu_h [[ \nabla_p (\vr_h s_h) ]] \ [[ p_h ]] } \dt\\
&= - c_v \int_0^T \sum_{ \sigma \in \Sigma_{int} } \intSh{ \mu_h \Ov{\vrh} [[ \vt_h ]]  \left[\left[ \frac{1}{\vt_h} \right]\right] } \dt
+ \int_0^T \sum_{ \sigma \in \Sigma_{int} } \intSh{ \mu_h [[\vrh]]\ [[\log(\vrh) ]] } \dt\\
&- c_v \int_0^T \sum_{ \sigma \in \Sigma_{int} } \intSh{ \mu_h \left(
[[\log(\vt_h) ]] + \Ov{\vt_h} \left[\left[ \frac{1}{\vt_h} \right]\right] \right) [[ \vrh]] } \dt.
\end{split}
\]
Next, we show that
\begin{equation} \label{nk6}
- [[ \vr_h ]] \left( [[ \log(\vt_h ) ]] + \Ov{\vt_h} \left[\left[ \frac{1}{\vt_h} \right]\right]  \right)
\leq - \frac{1}{2} |\ [[ \vr_h ]] \ | \ [[ \vt_h ]] \left[\left[ \frac{1}{\vt_h} \right]\right].
\end{equation}
As both expression in the above inequality are invariant with respect to the change  ``in'' and ``out'' and, in addition, the right--hand side is invariant with respect to the same operation in $\vrh$ and $\vt_h$ separately, it is enough to show (\ref{nk6}) assuming
$\vr_h^{\rm in} \geq \vr^{\rm out}_h$. In other words,
\[
- [[ \vr_h ]] = |\ [[ \vr_h ]] \ | \geq 0.
\]
Consequently, the proof of (\ref{nk6}) reduces to the inequality
\[
[[ \log(\vt_h ) ]] + \Ov{\vt_h} \left[\left[ \frac{1}{\vt_h} \right]\right] \leq - \frac{1}{2} [[ \vt_h ]] \left[\left[ \frac{1}{\vt_h} \right]\right].
\]
Denoting $Z = \frac{ \vt^{\rm out}_h }{\vt^{\rm in}_h }$, we have to show
\[
\log(Z) - \frac{1}{2} \left( Z - \frac{1}{Z} \right) \leq \frac{1}{2} \left( Z + \frac{1}{Z} \right) - 1
\]
or
\[
\log(Z) \leq Z - 1,
\]
which is obvious as $\log$ is a concave function.
In view of (\ref{nk6}), relation (\ref{as11}) yields
\begin{equation} \label{as12}
\begin{split}
h^{\alpha - 1}& \int_0^T \sum_{ \sigma \in \Sigma_{int} } \intSh{ [[ \vu_h ]]^2 \Ov{ \left( \frac{1}{\vt_h} \right) }  } \dt
+ \frac{1}{2} \int_0^T \sum_{ \sigma \in \Sigma_{int} } \intSh{  \left( \frac{1}{\vt_h} \right)^{\rm down}  \vr_h^{\rm up} |\Ov{\vu_h} \cdot \vc{n} |
[[ \vu_h ]]^2   } \dt \\
+ &\int_0^T \sum_{ \sigma \in \Sigma_{int} } \intSh{ \mu_h \Ov{\vr_h} [[\vu_h ]]^2 \Ov{ \left( \frac{1}{\vt_h} \right) } } \dt
+ \int_0^T \sum_{ \sigma \in \Sigma_{int} } \intSh{  \left( \widetilde{[[ b(\vrh) ]]} -
b'( {\vr_h^{\rm down} } ) \widetilde{[[ \vrh ]]} \right) | \Ov{\vuh} \cdot \vc{n} |} \dt
 \\
-  c_v &\int_0^T \sum_{ \sigma \in \Sigma_{int} } \intSh{    \vr_h^{\rm up}
\left( \widetilde{[[\log(\vt_h)]]} - \frac{1}{\vt_h^{\rm down}} \widetilde{[[ \vt_h ]]} \right) |\Ov{\vu_h} \cdot \vc{n} | } \dt \\
- c_v &\int_0^T \sum_{ \sigma \in \Sigma_{int} } \intSh{ \mu_h \min \{ \vr^{\rm in}_h, \vr^{\rm out}_h \} [[ \vt_h ]] \left[\left[ \frac{1}{\vt_h} \right]\right] } \dt
+ \int_0^T \sum_{ \sigma \in \Sigma_{int} } \intSh{ \mu_h [[\vrh]]\ [[\log(\vrh) ]] } \dt \\
&\aleq (1 + E_0).
\end{split}
\end{equation}

\section{Consistency}
\label{C}

We show consistency of the scheme (\ref{n2}--\ref{n4}), meaning the approximate solutions satisfy the weak formulation of the problem modulo approximation errors vanishing in the asymptotic limit $h \to 0$.

\subsection{Numerical flux}\label{CS1}

Firstly, we handle the numerical fluxes in \eqref{n2}, \eqref{n3} and the numerical entropy flux \eqref{num_flux_entropy} consisting of the upwind and $\mu_h$--dependent terms.

\subsubsection{Upwinds}

The upwind terms in the continuity equation (\ref{n2}), momentum equation (\ref{n3}), and the renormalized entropy balance (\ref{k6}) read
\begin{align*}
&\intOh{ \left( \Ov{\vr_h} \ \Ov{\vu_h} \cdot \vc{n} - \frac{1}{2} | \Ov{\vu_h} \cdot \vc{n} | [[\vr]]  \right) [[\Phi]]} ,\\
&\intOh{ \left( \Ov{\vc{m}_h}\ ( \Ov{\vu_h} \cdot \vc{n} ) - \frac{1}{2} | \Ov{\vu_h} \cdot \vc{n} | [[\vc{m}_h]] \right)  \cdot [[ \Phi ]] },\\
\mbox{and}\ &\intOh{ \left( \Ov{\vr_h \chi (s_h) } \Ov{\vu_h} \cdot \vc{n} - \frac{1}{2} | \Ov{\vu_h} \cdot \vc{n} | [[
\vr_h \chi(s_h) ]] \right) [[ \Phi ]] }, \ \mbox{respectively.}
\end{align*}
For $\Phi \in C^1(\Ov{\Omega_h})$ we get
\begin{align*}
&\intOh{ \vr_h b_h \vu_h \cdot \Grad \Phi } = \sum_{K \in \grid} \int_{K} \vr_h b_h \vu_h \cdot \Grad \Phi \ \dx =
\sum_{K \in \grid} \int_{\partial K} \vr_h b_h \vu_h \cdot \vc{n} \Phi \ {\rm d}S_h  \\
&= \sum_{K \in \grid} \int_{\partial K} \vr_h b_h \vu_h \cdot \vc{n} \left( \Phi - \Ov{\Pi_h [\Phi]}  \right) \ {\rm d}S_h
+ \sum_{K \in \grid} \int_{\partial K} \vr_h b_h \vu_h \cdot \vc{n} \Ov{\Pi_h [\Phi]}  \  {\rm d}S_h \\
&= - \sum_{\sigma \in \Sigma_{int}} \intSh{ [[ \vr_h b_h \vu_h ]] \cdot \vc{n} \left( \Phi - \Ov{\Pi_h [\Phi]}  \right) }
- \sum_{\sigma \in \Sigma_{int}} \intSh{ [[ \vr_h b_h \vu_h ]] \cdot \vc{n} \Ov{\Pi_h [\Phi] } }\\
&= - \sum_{\sigma \in \Sigma_{int}} \intSh{ [[ \vr_h b_h \vu_h ]] \cdot \vc{n} \left( \Phi - \Ov{\Pi_h [\Phi]}  \right) }
+ \sum_{\sigma \in \Sigma_{int}} \intSh{ \Ov{ \vr_h b_h \vu_h }  \cdot \vc{n} [[ \Pi_h [\Phi] ]] }\\
&= \sum_{\sigma \in \Sigma_{int}} \intSh{ \Ov{ \vr_h b_h } \ \Ov{\vu_h}  \cdot \vc{n} [[ \Pi_h [\Phi] ]] }
\\
&+ \sum_{\sigma \in \Sigma_{int}} \intSh{ \left( \Ov{ \vr_h b_h \vu_h}   - \Ov{ \vr_h b_h } \ \Ov{\vu_h} \right)
 \cdot \vc{n} [[ \Pi_h [\Phi] ]] }
- \sum_{\sigma \in \Sigma_{int}} \intSh{ [[ \vr_h b_h \vu_h ]] \cdot \vc{n} \left( \Phi - \Ov{\Pi_h [\Phi]}  \right) }\\
&= \sum_{\sigma \in \Sigma_{int}} \intSh{ Up[ \vr_h b_h ] \ [[\ \Pi[\Phi]\ ]] } + \frac{1}{2} \sum_{\sigma \in \Sigma_{int}} \intSh{
|\Ov{\vu_h} \cdot \vc{n} | [[\vr_h b_h ]] [[\ \Pi [\Phi] \ ]] } \\
&+ \sum_{\sigma \in \Sigma_{int}} \intSh{ \left( \Ov{ \vr_h b_h \vu_h}  - \Ov{ \vr_h b_h } \ \Ov{\vu_h} \right)
 \cdot \vc{n} [[ \Pi_h [\Phi] ]] }
- \sum_{\sigma \in \Sigma_{int}} \intSh{ [[ \vr_h b_h \vu_h ]] \cdot \vc{n} \left( \Phi - \Ov{\Pi_h [\Phi]}  \right) }.
\end{align*}
Seeing that
\[
\Ov{u v} - \Ov{u} \ \Ov{v} = \frac{1}{4} [[u]]\ [[v]]
\]
we have to control the following error terms:
\[
\begin{split}
E_1 &= \sum_{\sigma \in \Sigma_{int}} \intSh{
|\Ov{\vu_h} \cdot \vc{n} | [[\vr_h b_h ]] [[\ \Pi [\Phi] \ ]] },\\
E_2 &= \sum_{\sigma \in \Sigma_{int}} \intSh{ [[ \vr_h b_h \vu_h ]] \cdot \vc{n} \left( \Phi - \Ov{\Pi_h [\Phi]}  \right) },\\
E_3 &= \sum_{\sigma \in \Sigma_{int}} \intSh{
[[\vr_h b_h ]]\ [[\vu_h]] \cdot \vc{n}  [[\ \Pi [\Phi] \ ]] },
\end{split}
\]
where  $b_h$ is either 1 or $\chi(s_h)$ or $u^j_h,$ $ j=1, \dots, N$.
In view of (\ref{n4c}) and the identity
\[
[[ \vr_h b_h \vu_h ]]\cdot \vc{n} = [[ \vr_h b_h ]] \Ov{\vu_h} \cdot \vc{n} + \Ov{ \vr_h b_h } [[\vu_h]]\cdot \vc{n},
\]
it is enough to show that
\begin{equation} \label{C1-}
\begin{split}
E_{1,h}&= h \sum_{\sigma \in \Sigma_{int}} \intSh{
|\Ov{\vu_h} \cdot \vc{n} | \ \left|\ [[\vr_h b_h ]] \ \right|  } \to 0,\\
E_{2,h} &= h \sum_{\sigma \in \Sigma_{int}} \intSh{ | \Ov{\vr_h b_h} | \ | \ [[ \vu_h  ]]\cdot \vc{n}  \ |  }
\to 0,\\
E_{3,h} &= h \sum_{\sigma \in \Sigma_{int}} \intSh{
\left| \ [[\vr_h b_h ]]\ \right| \ \left| [[\vu_h]] \cdot \vc{n} \right|  } \to 0
\end{split}
\end{equation}
as $h \to 0$
for any fixed $\Phi \in C^1(\Ov{\Omega}_h)$. Moreover, by virtue of the minimum entropy principle (\ref{as4}), the entropy $s_h$ is bounded below uniformly for $h \to 0$. As the cut--off function $\chi$ is supposed to be bounded from above, we may assume
\begin{align*}
| \chi(s_h) | \aleq 1 \ \mbox{for}\ h \to 0.
\end{align*}

The following analysis leans heavily on the bound
\begin{equation} \label{C1}
\int_0^T \sum_{\sigma \in \Sigma_{int}} \intSh{ [[ \vu_h ]]^2 } \dt \aleq h^{1 - \alpha}
\end{equation}
that follows directly from the entropy estimates (\ref{as12}) provided
\begin{equation} \label{C2}
0 < \vt_h \aleq 1 \ \mbox{uniformly for}\ h \to 0.
\end{equation}
Accordingly, we \emph{suppose} that the approximate solutions satisfy (\ref{C2}). Then, as $\gamma > 1$, the entropy
minimum principle (\ref{as4}) yields a similar bound on the density,
\begin{equation} \label{C3}
0< \vr_h \aleq 1 \ \mbox{uniformly for}\ h \to 0.
\end{equation}

With (\ref{C2}), (\ref{C3}) at hand, the convergence of the errors $E_{2,h}$, $E_{3,h}$ for  $b_h = 1$ and $b_h=\chi(s_h)$ reduces to showing
\[
h \sum_{\sigma \in \Sigma_{int}} \intSh{ | \ [[ \vu_h ]] \ | } \to 0.
\]
To see this, we use H\" older's inequality,
\[
\begin{split}
h \sum_{\sigma \in \Sigma_{int}} &\intSh{ | \ [[ \vu_h ]] \ | }
\leq h \left( \sum_{\sigma \in \Sigma_{int}} \intSh{ [[ \vu_h ]]^2 } \right)^{1/2} \left( \sum_{\sigma \in \Sigma_{int}} \intSh{ 1 } \right)^{1/2}
\\
&\aleq \sqrt{h} \left( \sum_{\sigma \in \Sigma_{int}} \intSh{ [[ \vu_h ]]^2 } \right)^{1/2}
\aleq h^{1 - \frac{\alpha}{ 2}} F^1_h  ,\ \| F^1_h \|_{L^2(0,T)} \aleq 1,
\end{split}
\]
where the last inequality follows from the hypothesis (\ref{C1}).

In order to control the integral in $E_{1,h}$ we need bounds on the velocity $\vu_h$. They can be deduced from the total energy balance
(\ref{as5}) if we make another extra {\em hypothesis}, namely
\begin{equation} \label{C6}
0 < \underline{\vr} \leq \vr_h \ \mbox{uniformly for all}\ h \to 0.
\end{equation}
In view of (\ref{as4}) this implies a similar lower bound on the approximate temperature,
\begin{equation} \label{C7}
0 < \underline{\vt} \leq \vt_h \ \mbox{uniformly for all}\ h \to 0.
\end{equation}
Under these circumstances, we easily deduce from (\ref{as5}), (\ref{as12}) the following bounds:
\begin{equation} \label{C8}
\sup_{t \in [0,T]} \| \vu_h(t) \|_{L^2(\Omega_h)} \aleq 1,
\end{equation}
\begin{equation} \label{C9}
\int_0^T \sum_{\sigma \in \Sigma_{int}} \intSh{ |\Ov{\vu_h} \cdot \vc{n}| [[ \vu_h ]]^2 } \dt \aleq 1,
\end{equation}
\begin{equation} \label{C10}
\int_0^T \sum_{\sigma \in \Sigma_{int}} \intSh{ |\Ov{\vu_h} \cdot \vc{n}| [[ \vr_h ]]^2 } \dt \aleq 1,
\end{equation}
\begin{equation} \label{C11}
\int_0^T \sum_{\sigma \in \Sigma_{int}} \intSh{ |\Ov{\vu_h} \cdot \vc{n}| [[ \vt_h ]]^2 } \dt \aleq 1.
\end{equation}
In particular,  we obtain the estimates
\begin{subequations}\label{weakBV}
\begin{align}
h^{\alpha-1}&\int_0^T\sum_{ \sigma \in \Sigma_{int}} \intSh{[[\vu_h]]^2}\ \dt\aleq 1, \\
 \int_0^T\sum_{ \sigma \in \Sigma_{int}} \intSh{ \lambda_h [[\vrh]]^2 } \ \dt \aleq 1, \ &\int_0^T\sum_{ \sigma \in \Sigma_{int}} \intSh{ \lambda_h [[\vt_h]]^2 } \ \dt \aleq1, \ \lambda_h \approx |\Ov{\vu_h}\cdot\vn|+\mu_h,\label{weakBV_vr_vt}
\end{align}
\end{subequations}
which are slightly better than the standard weak BV estimates, cf. \cite{FLM18,FjMiTa2,FjMiTa1}.

Now, the error term $E_{1,h}$ for  $b_h$ either equal to 1 or $\chi(s_h)$ can be handled as
\[
\begin{split}
h \sum_{\sigma \in \Sigma_{int}} &\intSh{
|\Ov{\vu_h} \cdot \vc{n} | \ \left|\ [[\vr_h b_h ]] \ \right|  } \aleq
h \sum_{\sigma \in \Sigma_{int}} \intSh{
|\Ov{\vu_h} \cdot \vc{n} | \ |\ [[\vr_h ]] \ | + | \ [[\vt_h ]] \ |  }\\
&\aleq h \left( \sum_{\sigma \in \Sigma_{int}} \intSh{
|\Ov{\vu_h}| } \right)^{1/2}\left( \sum_{\sigma \in \Sigma_{int}} \intSh{
|\Ov{\vu_h} \cdot \vc{n} |\left( [[\vr_h ]]^2 +  [[\vt_h ]]^2 \right) } \right)^{1/2} \\
&\aleq \sqrt{h} \| \vu_h \|_{L^1(\Omega_h)}^{1/2} F^2_h\aleq \sqrt{h} F^2_h, \  \| F_h^2 \|_{L^2(0,T)} \aleq 1.
\end{split}
\]
Thus it remains to estimate $E_{1,h}$, $E_{2,h}$, $E_{3,h}$ for $b_h = u^j_h$. For $E_{2,h}$, we get
\begin{equation} \label{C13}
\begin{split}
h \sum_{\sigma \in \Sigma_{int}} & \intSh{ | \Ov{\vr_h u^j_h} | \ | \ [[ \vu_h  \cdot \vc{n} ]] \ |  } \aleq
h \left( \sum_{\sigma \in \Sigma_{int}} \intSh{ [[ \vu_h ]]^2 } \right)^{1/2}
\left( \sum_{\sigma \in \Sigma_{int}} \intSh{ |\Ov{\vu}_h |^2 } \right)^{1/2} \\
&\aleq h^{\frac{3}{2} - \frac{\alpha}{2}}
 F^1_h\, h^{-\frac{1}{2}}\|\vu_h\|_{L^2(\Omega_h)}
\leq h^{1 - \frac{\alpha}{2}} F^1_h, \ \| F^1_h \|_{L^2(0,T)} \aleq 1,
\end{split}
\end{equation}
where we have used the trace inequality, (\ref{C1}) and (\ref{C8}). As for $E_{3,h}$ it rewrites as
\[
\begin{split}
h \sum_{\sigma \in \Sigma_{int}} & \intSh{
\left| \ [[\vr_h u^j_h ]]\ \right| \ \left| [[\vu_h]] \cdot \vc{n} \right|  }
\aleq h \sum_{\sigma \in \Sigma_{int}} \intSh{
[[ \vu_h ]]^2  } + h \sum_{\sigma \in \Sigma_{int}} \intSh{
 |\Ov{\vu_h}|\ |\ [[ \vu_h ]] \ |  }\\
&\aleq h^{2 - \alpha} F^3_h + h \sum_{\sigma \in \Sigma_{int}} \intSh{
 |\Ov{\vu_h}|\ |\ [[ \vu_h ]] \ |  },\ \|F^3_h \|_{L^1(0,T)} \aleq 1,
\end{split}
\]
while the last integral can be handled exactly as in (\ref{C13}).
Finally, we are left with $E_{1,h}$, specifically,
\[
\begin{split}
h \sum_{\sigma \in \Sigma_{int}} &\intSh{
|\Ov{\vu_h} \cdot \vc{n} | \ \left|\ [[\vr_h u^j_h ]] \ \right|  }
\\ &\aleq h \sum_{\sigma \in \Sigma_{int}} \intSh{
|\Ov{\vu_h} \cdot \vc{n}| \ |\Ov{\vu}_h | \ |\ [[\vr_h ]] \ |   }
+h \sum_{\sigma \in \Sigma_{int}} \intSh{
|\Ov{\vu_h} \cdot \vc{n}| |\ [[\vu_h ]] \ |   },
\end{split}
\]
where the last integral can be estimated exactly as in (\ref{C13}). Next, by H\" older's inequality, the trace inequality,  (\ref{C10}), we get
\[
\begin{split}
h \sum_{\sigma \in \Sigma_{int}} &\intSh{
|\Ov{\vu_h} \cdot \vc{n}| \ |\Ov{\vu}_h | \ |\ [[\vr_h ]] \ |   }
\aleq h \left( \sum_{\sigma \in \Sigma_{int}} \intSh{ | \Ov{\vc{u}_h} |^3  } \right)^{1/2}\left(\sum_{\sigma \in \Sigma_{int}} \intSh{ | \Ov{\vu_h} \cdot \vc{n}| [[ \vr_h ]]^2 } \right)^{1/2}
\\
&\aleq \sqrt{h} \| \vu_h \|_{L^3(\Omega_h)}^{3/2} \,  F^4_h, \ \| F^4_h \|_{L^2(0,T)} \aleq 1.
\end{split}
\]
Now, in view of the interpolation inequality
\[
\| \vu_h \|_{L^3(\Omega_h)} \aleq \| \vu_h \|_{L^2(\Omega_h)}^{1/2} \| \vu_h \|_{L^6(\Omega_h)}^{1/2},
\]
combined with (\ref{C8}),
we obtain
\[
h \sum_{\sigma \in \Sigma_{int}} \intSh{
|\Ov{\vu_h} \cdot \vc{n}| \ |\Ov{\vu}_h | \ |\ [[\vr_h ]] \ |   }
\aleq \sqrt{h}  \| \vu_h \|_{L^6(\Omega_h)}^{3/4}\, F^4_h.
\]
Finally, we apply the discrete Sobolev embedding (\ref{n4d}) and (\ref{C1}) to conclude
\[
\begin{split}
h \sum_{\sigma \in \Sigma_{int}} \intSh{
|\Ov{\vu_h} \cdot \vc{n}| \ |\Ov{\vu}_h | \ |\ [[\vr_h ]] \ |   }
&\aleq \sqrt{h} F^4_h \left( 1 +  \left( \sum_{\sigma \in \Sigma_{int}} \intSh{
\frac{[[ \vu_h ]]^2}{h} } \right)^{1/2} \right)^{3/4}\aleq h^{ \frac{4 - 3 \alpha}{8} } F^4_h F^5_h,
\end{split}
\]
with
\[
\| F^5_h \|_{L^{8/3}(0,T)} \aleq 1.
\]
Thus the error in the upwind terms satisfies (\ref{C1-}) as soon as
\[
0 < \alpha < \frac{4}{3},
\]
and the extra hypotheses (\ref{C2}), (\ref{C6}) hold.

\subsubsection{$\mu_h$--dependent terms}

The numerical fluxes of the continuity and momentum  equations \eqref{n2}, \eqref{n3}, and the numerical entropy flux \eqref{num_flux_entropy} contain  $\mu_h$--dependent terms, namely 
\begin{align*}
&\sum_{\sigma \in \Sigma_{int}} \intSh{\mu_h [[\vr_h]] [[\ \Pi[\Phi]\ ]] }, \\
&\sum_{\sigma \in \Sigma_{int}} \intSh{\mu_h   [[\vr_h\vu_h]] \cdot [[\ \Pi[\bfPhi]\ ]]}, \\
&\sum_{\sigma \in \Sigma_{int}} \intSh{\mu_h \bigg( \Ov{\nabla_{\vr}(\vr_h\chi(s_h))}  \ [[\vr_h]]  + \Ov{\nabla_{p}(\vr_h\chi(s_h))}\ [[p_h]]  \bigg) [[\ \Pi[\Phi]\ ]]}.
\end{align*}
In what follows we show  they vanish in the limit $h \rightarrow 0.$ 
In view of our hypotheses (\ref{C2}), (\ref{C6}), the product rule yields 
\begin{align}\label{p_vr_vt}
[[ p_h ]] \approx [[\vr_h ]] + [[ \vt_h ]],
\end{align}
and the estimates \eqref{weakBV_vr_vt} imply
\begin{align}
\int_0^T\sum_{\sigma\in\Sigma_{int}}\intSh{\mu_h[[\vr_h]]^2} \aleq 1,\label{C14} \\
\int_0^T\sum_{\sigma\in\Sigma_{int}}\intSh{\mu_h[[\vt_h]]^2} \aleq 1.\label{C15}
\end{align}

{\em Assuming}  the parameter  $\mu_h$ is bounded, H\" older's inequality with \eqref{n4c} and \eqref{C14} directly yield
\begin{align*}
\sum_{\sigma \in \Sigma_{int}} &\intSh{\mu_h[[\vr_h]]\ [[\ \Pi [\Phi]\ ]]} \aleq \left(\sum_{\sigma \in \Sigma_{int}} \intSh{\mu_h [[\vr_h]]^2  }\right)^{1/2}\left(\sum_{\sigma\in\Sigma_{int}}\intSh{\mu_h h^2} \right)^{1/2} \\
&\aleq \sqrt{h}F_h^6, \ \|F_h^6\|_{L^2(0,T)} \aleq 1.
\end{align*}
Analogously, using the product rule, the trace inequality, and  bounds \eqref{C1}, \eqref{C3}, \eqref{C14}, we get
\begin{align*}
\sum_{\sigma \in \Sigma_{int}} &\intSh{\mu_h[[\vr_h\vu_h]]\cdot [[\ \Pi [\bfPhi]\ ]]}\\
& \aleq h\left(\sum_{\sigma\in\Sigma_{int}}\intSh{\mu_h[[\vr_h]]^2}\right)^{1/2} \left(\sum_{\sigma \in \Sigma_{int}} \intSh{\mu_h|\Ov{\vu_h}|^2}\right)^{1/2} +\sqrt{h}\left( \sum_{\sigma \in \Sigma_{int}} \intSh{[[\vu_h]]^2}\right)^{1/2} \\
& \aleq \sqrt{h}\|\vu_h\|_{L^2(\Omega_h)}F_h^6 + h^{1 - \frac{\alpha}{ 2}} F^1_h \aleq \sqrt{h}F_h^6 +h^{1-\frac{\alpha}{2}}F_h^1, \ \|F_h^1\|_{L^2(0,T)},\ \|F_h^6\|_{L^2(0,T)} \aleq 1.
\end{align*}
Finally,  \eqref{p_vr_vt}, \eqref{C14} and \eqref{C15}  imply
\begin{align*}
\sum_{ \sigma \in \Sigma_{int} } &\intSh{ \mu_h \left( \Ov{\nabla_{\vr} (\vr_h \chi(s_h))}[[\vr_h ]]+\Ov{\nabla_p (\vr_h \chi(s_h))}[[p_h ]]\right)\ [[\ \Pi [\Phi]\ ]]  } \\
& \aleq \sqrt{h}\left(\sum_{\sigma \in \Sigma_{int}}\intSh{\mu_h[[\vr_h]]^2}\right)^{1/2}+\sqrt{h}\left(\sum_{\sigma \in \Sigma_{int}}\intSh{\mu_h[[\vt_h]]^2}\right)^{1/2} \\
&\aleq \sqrt{h}(F_h^6 + F_h^7),\ \|F_h^6\|_{L^2(0,T)},\ \|F_h^7\|_{L^2(0,T)} \aleq 1.
\end{align*}

\subsection{The artificial viscosity and the pressure terms}
\label{CS2}

There are two remaining terms to be handled in the momentum equation, namely,
\begin{align*}
&h^{\alpha - 1} \sum_{\sigma \in \Sigma_{int}} \intSh{ [[ \vu_h]] \cdot [[ \ \Pi [\bfPhi] \ ]] },
\end{align*}
and
\begin{align*}
&\sum_{\sigma \in \Sigma_{int}} \intSh{ \Ov{p}_h \vc{n} \cdot [[\ \Pi [ \bfPhi ] \ ]] }.
\end{align*}
First, in accordance with (\ref{n4c}),
\[
\begin{split}
 h^{\alpha - 1} & \left| \sum_{\sigma \in \Sigma_{int}} \intSh{ [[ \vu_h]] \cdot [[ \ \Pi [\bfPhi] \ ]] } \right|
\aleq h^{\alpha - 1} \left( \sum_{\sigma \in \Sigma_{int}} \intSh{ [[ \vu_h]]^2 } \right)^{1/2}
\left( \sum_{\sigma \in \Sigma_{int}} \intSh{ h^2 } \right)^{1/2} \\
&\aleq h^{\alpha - \frac{1}{2}} \left( \sum_{\sigma \in \Sigma_{int}} \intSh{ [[ \vu_h]]^2 } \right)^{1/2}
\aleq h^{\alpha/2} F^7_h, \ \| F^7_h \|_{L^2(0,T)} \aleq 1.
\end{split}
\]
Second,
\[
\begin{split}
&\intOh{ p_h \Div \bfPhi } = \sum_{K \in \grid} \int_K p_h \Div \bfPhi \ \dx =
\sum_{K \in \partial K_h} \int_{\partial K} p_h \bfPhi \cdot \vc{n} {\rm d}S_h \\
&= \sum_{\sigma \in \Sigma_{int}} \intSh{ [[p_h]] \left( \bfPhi  - \Ov{ \Pi[ \bfPhi ] } \right) \cdot \vc{n} }
+ \sum_{\sigma \in \Sigma_{int}} \intSh{ [[p_h]] \Ov{ \Pi[ \bfPhi ] }  \cdot \vc{n} }
\\
&=\sum_{\sigma \in \Sigma_{int}} \intSh{ [[p_h]] \left( \bfPhi  - \Ov{ \Pi[ \bfPhi ] } \right) \cdot \vc{n} }
- \sum_{\sigma \in \Sigma_{int}} \intSh{ \Ov{p_h } [[ \Pi[ \bfPhi ] ]]  \cdot \vc{n} } .
\end{split}
\]
Here, similarly to the preceding section, the error term can be estimated as
\[
\begin{split}
\left| \sum_{\sigma \in \Sigma_{int}} \intSh{ [[p_h]] \left( \bfPhi  - \Ov{ \Pi[ \bfPhi ] } \right) \cdot \vc{n} }
\right| &\aleq \left( \sum_{\sigma \in \Sigma_{int}} \intSh{ [[p_h]]^2 } \right)^{1/2} \left( \sum_{\sigma \in \Sigma_{int}}
\intSh{ h^2 } \right)^{1/2} \\ &\aleq
\sqrt{h} \left( \sum_{\sigma \in \Sigma_{int}} \intSh{ [[p_h]]^2 } \right)^{1/2}.
\end{split}
\]
Recall that $
[[ p_h ]] \approx [[\vr_h ]] + [[ \vt_h ]];$
whence for the error to tend to zero it is enough to {\em assume}
\[
\mu_h \ageq h^\beta >0, \ 0 \leq \beta < 1 .
\]
In  the case of uniform rectangular/cubic elements  we allow $\mu_h = 0$. Indeed, due to \eqref{n4c2} and \eqref{trace} we have, for any $\bfPhi \in C^2(\Ov{\Omega}_h;R^N),$
\begin{align*}
\left| \sum_{\sigma \in \Sigma_{int}} \intSh{ [[p_h]] \left( \bfPhi  - \Ov{ \Pi[ \bfPhi ] } \right) \cdot \vc{n} }\right| 
&\aleq \sum_{\sigma \in \Sigma_{int}} \left|\ [[p_h]]\ \right| \intSh{\left| \bfPhi  - \Ov{ \Pi[ \bfPhi ]} \right|}    \\
&\aleq \sum_{\sigma \in \Sigma_{int}}\left(\intSh{ |\ [[p_h]]\ | }\right)\ h^2 \aleq  h \sum_{K \in \grid} \int_{K} |p_h| {\rm d}x  ,
\end{align*}
which tends to zero as $p_h \in L^\infty((0,T); L^1(\Omega_h))$.

\subsection{Consistency formulation}

Summing up the  results of Subsections~\ref{CS1} and \ref{CS2}, we obtain a consistency formulation of the approximation scheme (\ref{n2}--\ref{n4}).

\begin{Theorem} \label{Tm2}

Let the initial data $\vr_{0,h}$, $\vc{m}_{0,h}$, $E_{0,h}$ satisfy the hypotheses of Theorem \ref{Tm1}. Let $[\vr_h, \vc{m}_h, E_h]$
be the unique solutions of the approximate problem (\ref{n2}--\ref{n4}) on the time interval $[0,T]$.

Then
\begin{equation} \label{cP1}
\left[ \intOh{ \vr_h \varphi } \right]_{t=0}^{t = \tau} =
\int_0^\tau \intOh{ \left[ \vr_h \partial_t \varphi + \vc{m}_h \cdot \Grad \varphi \right]} \dt  + \int_0^T
e_{1,h} (t, \varphi) \dt
\end{equation}
for any $\varphi \in C^1([0,T] \times \Ov{\Omega}_h)$;
\begin{equation} \label{cP2}
\left[ \intOh{ \vc{m}_h \bfphi } \right]_{t=0}^{t = \tau} =
\int_0^\tau \intOh{ \left[ \vc{m}_h \cdot \partial_t \bfphi + \frac{\vc{m}_h \otimes \vc{m}_h} {\vr_h} : \Grad \bfphi
+ p_h \Div \bfphi \right]} \dt  + \int_0^T
e_{2,h} (t, \bfphi) \dt
\end{equation}
for any $\bfphi \in C^1([0,T] \times \Ov{\Omega}_h; R^N)$, $\bfphi \cdot \vc{n}|_{\Omega_h} = 0$;
\begin{equation} \label{cP3}
\intOh{ E_h(t) } = \intOh{ E_{0,h} };
\end{equation}
\begin{equation} \label{cP4}
\left[ \intOh{ \vr_h \chi(s_h) \varphi } \right]_{t=0}^{t = \tau} \geq
\int_0^\tau \intOh{ \left[ \vr_h \chi (s_h) \partial_t \varphi + \chi(s_h) \vc{m}_h \cdot \Grad \varphi \right]} \dt  + \int_0^T
e_{3,h} (t, \varphi) \dt
\end{equation}
for any $\varphi \in C^1([0,T] \times \Ov{\Omega}_h)$, $\varphi \geq 0$, and any $\chi$,
\[
\chi : R \to R \ \mbox{a non--decreasing concave function}, \ \chi(s) \leq \Ov{\chi} \ \mbox{for all}\ s \in R.
\]

If, in addition,
\begin{equation} \label{cP5}
 h^\beta \aleq \mu_h \aleq 1, \ 0 \leq \beta < 1,\ 0< \alpha < \frac{4}{3},
\end{equation}
and
\begin{equation} \label{cP6}
0 < \Ov{\vr} \leq \vr_h(t), \ \vt_h(t) \leq \Ov{\vt} \ \mbox{for all}\ t \in [0,T] \ \mbox{uniformly for}\ h \to 0,
\end{equation}
then
\[
\| e_{j,h} (\cdot, \varphi ) \|_{L^1(0,T)} \aleq h^\delta \| \varphi \|_{C^1} \ \mbox{for some}\ \delta > 0.
\]
\end{Theorem}
 In the case of  uniform rectangular/cubic elements the result of Theorem~\ref{Tm2} remains valid  for $0\leq \mu_h\aleq 1,$ and $\bfphi \in C^1([0,T]; C^2(\Ov{\Omega}_h; R^N))$, $\bfphi \cdot \vc{n}|_{\Omega_h} = 0.$

\begin{Remark} \label{NR2}

Omitting the $h^{\alpha-1}$--dependent terms in (\ref{n3}--\ref{n4}) corresponds to the Lax--Friedrichs scheme  with the numerical fluxes
\[
\Ov{r_h \vu_h} \cdot \vc{n} -   \lambda_h  [[ r_h ]] = \Ov{r_h }\ \Ov{\vu_h} \cdot \vc{n} -   \lambda_h  [[ r_h ]]  +\frac{1}{4}[[r_h][[\vu_h]] \cdot \vn
\]
with $\lambda_h = \frac{1}{2}\max(\lambda_h^{\rm in}, \lambda_h^{\rm out})$, $\lambda_h = \frac{1}{2}|\Ov{\vu_h} \cdot \vc{n}| + c_h$, where  $c_h= \sqrt{\gamma \vt_h}$  stands for the speed of sound.
 In the standard Lax--Friedrichs scheme the average $\Ov{r_h\vu_h}$ instead of $\Ov{r_h}\ \Ov{\vu_h}$ is used.
 Moreover, in the energy equation $\Ov{p_h\vuh}$ is used instead of (\ref{pressure}) for the pressure term in the energy flux, cf. Remark~\ref{nR1a}.  
Nevertheless, the present proof under the hypotheses \eqref{cP5}, \eqref{cP6} might be adapted to the standard Lax--Friedrichs scheme.
\end{Remark}

\section{Convergence}
\label{L}

In view of the uniform bounds (\ref{as2}), (\ref{as5}), and (\ref{as6}), the family $\{ \vrh, \vc{m}_h, E_h \}_{h > 0}$  of approximate solutions
is uniformly bounded in $L^\infty(0,T; L^1(\Omega_h))$. Moreover,
$\{ \vrh \}_{h > 0}$ is bounded in $L^\infty(0,T; L^{\gamma}(\Omega_h))$ and
$\{ \vc{m}_h \}_{h > 0}$ is bounded in $L^\infty (0,T; L^{\frac{2\gamma}{\gamma + 1}}(\Omega_h; R^N))$ uniformly for $h \to 0$.

\subsection{Young measure generated by the approximate solutions}

In accordance with the fundamental theorem on Young measures, see Ball \cite{BALL2} or Pedregal \cite{PED1},
the family $\{ \vrh, \vc{m}_h, E_h \}_{h > 0}$, up to a suitable subsequence, generates a Young measure
$\{ \mathcal{V}_{t,x} \}_{(t,x) \in (0,T) \times \Omega_h}$. Recall that the Young measure is an object with the following properties:
\begin{itemize}
\item the mapping
\[
\mathcal{V}_{t,x}:
(t,x) \in (0,T) \times \Omega_h \mapsto \mathcal{P}(\mathcal{F})
\]
is weakly-(*) measurable, where $\mathcal{P}$ is the space of probability measures defined on the phase space
\[
\mathcal{F} = \left\{ \vr, \vc{m}, E \ \Big| \ \vr \geq 0, \ \vc{m} \in R^N,\ E \geq 0 \right\};
\]
\item
\[
G(\vr_h, \vc{m}_h, E_h ) \to \Ov{G (\vr, \vc{m}, E)} \ \mbox{weakly-(*) in}\ L^\infty((0,T) \times \Omega_h)
\]
for any $G \in C_c( \mathcal{F} )$, and
\[
\Ov{G (\vr, \vc{m}, E)} (t,x) = \int_{\mathcal{F}} G(\vr, \vc{m}, E) {\rm d}\mathcal{V}_{t,x} \equiv
\left< \mathcal{V}_{t,x}; G(\vr, \vc{m}, E) \right> \ \mbox{for a.a.} \ (t,x) \in (0,T) \times \Omega_h.
\]

\end{itemize}

We shall use the following result proved in \cite[Lemma~2.1]{FGSWW1}.
\begin{Lemma}\label{lemma}
Let 
\begin{align*}
|G(\vr,\vm,E)|\leq F(\vr,\vm,E) \ \mbox{ for all }\ (\vr,\vm,E) \in \mathcal{F}.
\end{align*}
Then
\begin{align*}
\left| \Ov{G(\vr,\vm,E)} - \left< \mathcal{V}_{t,x}; G(\vr, \vc{m}, E) \right> \right| \leq \Ov{F(\vr,\vm,E)} -  \left< \mathcal{V}_{t,x}; F(\vr, \vc{m}, E) \right> \equiv \mu_F \mbox{ in } \mathcal{M}([0,T]\times\Omega).
\end{align*}
\end{Lemma}

\subsection{Kinetic energy concentration defect}

Under the extra hypotheses (\ref{cP6}), the support of the measure $\mathcal{V}_{t,x}$ is contained in the set
\[
{\rm supp} [ \mathcal{V}_{t,x} ] \subset \left\{ [\vr, \vc{m}, E] \ \Big| \
0 < \underline{\vr} \leq \vr \leq \Ov{\vr},\ 0 < \underline{\vt} \leq \vt \leq \Ov{\vt} \right\}
\ \mbox{for a.a.}\ (t,x) \in (0,T) \times \Omega_h.
\]
In particular, all non--linearities appearing in the consistency formulation (\ref{cP1} -- \ref{cP4}) are weakly precompact in the
Lebesgue space $L^1((0,T) \times \Omega_h)$, with the only exception of the convective term
\[
\left\{ \frac{\vc{m}_h \otimes \vc{m}_h }{\vr_h} \right\}_{h > 0}
\ \mbox{bounded in}\ L^1((0,T) \times \Omega_h, R^{N \times N}).
\]
For the latter we can only assert that
\[
\frac{\vc{m}_h \otimes \vc{m}_h }{\vr_h} \to \Ov{\frac{\vc{m} \otimes \vc{m} }{\vr}} \ \mbox{weakly-(*) in}\ \mathcal{M}([0,T] \times \Ov{\Omega}_h; R^{N \times N}).
\]
We denote
\[
\mathbb{C}_d = \Ov{\frac{\vc{m} \otimes \vc{m} }{\vr}}-\left< \mathcal{V}_{t,x}; \frac{\vc{m} \otimes \vc{m} }{\vr} \right> \in \mathcal{M}([0,T] \times \Ov{\Omega}_h;
R^{N \times N})
\]
the associated \emph{concentration defect measure}. As
\[
\left| \frac{\vc{m}_h \otimes \vc{m}_h }{\vr_h} \right| \aleq \frac{|\vc{m}_h|^2}{\vr_h} \leq E_h,
\]
we may use Lemma~\ref{lemma}
 to conclude that
\begin{equation} \label{L1}
\int_0^\tau \int_{\Ov{\Omega}_h } 1 \ {\rm d}|\mathbb{C}_d| \aleq \int_{\Omega_h} E_0 \dx - \intOh{ \left< \mathcal{V}_{\tau, x};
E \right> } \ \mbox{for a.a.}\ \tau \in [0,T].
\end{equation}
The quantity on the right--hand side of (\ref{L1}) is called \emph{energy dissipation defect} and inequality (\ref{L1}) plays a crucial
role in the concept of \emph{dissipative measure--valued (DMV) solutions} to the complete Euler system introduced in \cite{BreFei17}.

\subsection{Limit problem}

We say that a family of probability measures $\{ \mathcal{V}_{t,x} \}_{(t,x) \in (0,T) \times \Omega_h}$ is a (DMV) solution to the
complete Euler system (\ref{E7}--\ref{E9}) if:

\begin{itemize}
\item
\[
\left[ \intOh{ \left< \mathcal{V}_{t,x};  \vr \right> \varphi } \right]_{t=0}^{t = \tau} =
\int_0^\tau \intOh{ \left[ \left< \mathcal{V}_{t,x} ; \vr \right> \partial_t \varphi +
\left< \mathcal{V}_{t,x}; \vc{m} \right> \cdot \Grad \varphi \right]} \dt
\]
for any $\varphi \in C^1([0,T] \times \Ov{\Omega}_h)$;
\item
\[
\begin{split}
&\left[ \intOh{ \left< \mathcal{V}_{t,x}; \vc{m} \right> \cdot \bfphi } \right]_{t=0}^{t = \tau} \\ &=
\int_0^\tau \intOh{ \left[ \left< \mathcal{V}_{t,x}; \vc{m} \right> \cdot \partial_t \bfphi +
\left< \mathcal{V}_{t,x}; \frac{\vc{m} \otimes \vc{m}} {\vr} \right> : \Grad \bfphi
+ \left< \mathcal{V}_{t,x}; p \right> \Div \bfphi \right]} \dt \\
&+\int_0^\tau \int_{\Ov{\Omega}_h} \Grad \varphi : {\rm d}\mathbb{C}_d
\end{split}
\]
for any $\bfphi \in C^1([0,T] \times \Ov{\Omega}_h; R^N)$, $\bfphi \cdot \vc{n}|_{\Omega_h} = 0$;
\item
\[
\intOh{ \left< \mathcal{V}_{\tau,x}; E \right> } \leq \intOh{ E_{0} }
\]
for a.a. $\tau \in [0,T]$;
\item
\[
\left[ \intOh{ \left< \mathcal{V}_{t,x}; \vr \chi(s) \right> \varphi } \right]_{t=0}^{t = \tau} \geq
\int_0^\tau \intOh{ \left[ \left< \mathcal{V}_{t,x}; \vr \chi (s) \right> \partial_t \varphi + \left< \mathcal{V}_{t,x}; \chi(s) \vc{m}
\right> \cdot \Grad \varphi \right]} \dt
\]
for any $\varphi \in C^1([0,T] \times \Ov{\Omega}_h)$, $\varphi \geq 0$, and any $\chi$,
\[
\chi : R \to R \ \mbox{a non--decreasing concave function}, \ \chi(s) \leq \Ov{\chi} \ \mbox{for all}\ s \in R;
\]
\item
\[
\int_0^\tau \int_{\Ov{\Omega}_h } 1 \ {\rm d}|\mathbb{C}_d| \aleq \int_{\Omega_h} E_0 \dx - \intOh{ \left< \mathcal{V}_{\tau, x};
E \right> }
\]
for a.a. $\tau \in [0,T]$.
\end{itemize}

Summing up the preceding discussion, we can state the following result.

\begin{Theorem} \label{Tm3}
Let the initial data $\vr_{0,h}$, $\vm_{0,h}$, $E_{0,h}$ satisfy
\[
\vr_{0,h} \geq \underline{\vr} > 0,\ E_{0,h} - \frac{1}{2} \frac{ |\vc{m}_{0,h}|^2 }{\vr_{0,h}} > 0.
\]
Let  $[ \vr_h, \vc{m}_h, E_h ]$ be the solution of the scheme (\ref{n2}--\ref{n4}) such that
\begin{equation*}
 h^\beta \aleq \mu_h \aleq 1, \ 0 \leq \beta < 1,\ 0< \alpha < \frac{4}{3},
\end{equation*}
and
\begin{align*}
0 < \Ov{\vr} \leq \vr_h(t), \ \vt_h(t) \leq \Ov{\vt} \ \mbox{for all}\ t \in [0,T] \ \mbox{uniformly for}\ h \to 0.
\end{align*}
Then the family of approximate solutions  $\{ \vr_h, \vc{m}_h, E_h \}_{h > 0}$ generates a Young
measure $\{ \mathcal{V}_{t,x} \}_{(t,x) \in (0,T) \times \Omega_h}$ that is a (DMV) solution of the complete Euler system (\ref{E7}--\ref{E9}).
\end{Theorem}

Finally, evoking the weak (DMV)--strong uniqueness result proved in \cite[Theorem 3.3]{BreFei17} we conclude with the following corollary.

\begin{Corollary} \label{Tm4}

In addition to the hypotheses of Theorem \ref{Tm3}, suppose that the complete Euler system (\ref{E7}--\ref{E9}) admits a Lipschitz--continuous
solution $[\vr, \vc{m}, E]$ defined on $[0,T]$.

Then
\[
\vr_h \to \vr,\ \vc{m}_h \to \vc{m}, \ E_h \to E
\ \mbox{(strongly) in}\ L^1((0,T) \times \Omega_h).
\]

\end{Corollary}

\section*{Conclusion}
In the present paper we have studied the convergence of a new finite volume method for  multi--dimensional Euler equations
of gas dynamics. As the Euler system admits highly oscillatory solutions, in particular they are ill--posed in the class of weak entropy solutions for $L^\infty$--initial data \cite{FKKM17}, it is more natural to investigate the convergence in the class of  {\sl dissipative measure--valued (DMV) solutions}.  The (DMV) solutions represent the most general class of solutions that still satisfy the weak--strong uniqueness property. Thus, if the strong solution exists the (DMV) solution coincides with the strong one on its lifespan, cf.~\cite{BreFei17}.

Our study is inspired by the work of Guermond and Popov \cite{GuePop} who proposed a viscous regularization of the compressible
Euler equations satisfying the minimum entropy principle and positivity preserving properties. They also showed the connection to the two--velocities Brenner's model \cite{BREN2, BREN, BREN1}, which is a base of our new finite volume method
(\ref{n2}--\ref{n4}). The method is (i) positivity preserving, i.e.  discrete density, pressure and temperature are positive on any finite time interval, (ii) entropy stable and (iii) satisfies the
minimum entropy principle. Moreover, the discrete entropy inequality allows us to control certain weak BV--norms, cf.~(\ref{weakBV}).
These results together with a priori estimates (\ref{as1}--\ref{as9}) yield the consistency of the new finite volume method under mild hypothesis. 
Indeed, instead of conventional convergence results based on rather unrealistic hypothesis on uniform
boundedness of all physical quantities, we only require that the discrete temperature is bounded and vacuum does not appear, cf.~(\ref{cP6}). 
In Theorem~6.1 we have shown that the numerical solutions of the finite volume method (\ref{n2}--\ref{n4}) generate the (DMV) solution of
the Euler equations. Consequently, using the recent result on the (DMV)--strong uniqueness, we have proven the convergence to the strong solution on its lifespan. \\

It seems that the hypothesis on $\vr_h$ can be relaxed, though removing the
boundedness of $\vt_h$ remains open. This can be an interesting question for future study. Moreover, in order to preserve
the Galilean invariance of the Brenner model (\ref{b4}--\ref{b6}) it is possible to consider the symmetric gradient in the $h^\alpha$--diffusion terms and the same convergence result can be shown. As far as we know the present convergence result is the first result in the literature, where the convergence of a finite volume method has been proven for multi--dimensional Euler equations assuming only that the gas remains in its non--degenerate region.

\def\cprime{$'$} \def\ocirc#1{\ifmmode\setbox0=\hbox{$#1$}\dimen0=\ht0
  \advance\dimen0 by1pt\rlap{\hbox to\wd0{\hss\raise\dimen0
  \hbox{\hskip.2em$\scriptscriptstyle\circ$}\hss}}#1\else {\accent"17 #1}\fi}

\end{document}